\documentstyle[a4,12pt,amsfonts,leqno]{article}

\newtheorem{Theorem}{Theorem}[section]
\newtheorem{Proposition}[Theorem]{Proposition}
\newtheorem{Lemma}[Theorem]{Lemma}
\newtheorem{Corollary}[Theorem]{Corollary}
\newtheorem{Remark}[Theorem]{Remark}
\newtheorem{Example}[Theorem]{Example}
\newtheorem{Definition}[Theorem]{Definition}

\newtheorem{Claim}[Theorem]{Claim}

\newcounter{ctr}[section]

\newcommand{\bangou}{\addtocounter{ctr}{1}}

\def\RMN#1{\uppercase\expandafter{\romannumeral#1}}

\newcommand{\ya}[1]{\mathrel{\overrightarrow{\hphantom{#1}}}}
\newcommand{\uetuki}[2]{{\mathop{\scriptstyle {#1}}\limits_
{\scriptstyle {#2}}}}

\newcommand{\qed}{{\unskip\nobreak\hfil\penalty50\quad\null\nobreak\hfil{\bf
q.e.d.}\parfillskip0pt\finalhyphendemerits0\par\medskip}}
\newcommand{\proof}{\noindent{\it Proof.} \ }

\newcommand{\tor}[4]{{\rm Tor}_{#1}^{#2}({#3},{#4})}

\newcommand{\supp}{\mathop{\rm Supp}\nolimits}
\newcommand{\rank}{\mathop{\rm rank}\nolimits}

\newcommand{\spec}{\mathop{\rm Spec}\nolimits}
\newcommand{\proj}{\mathop{\rm Proj}\nolimits}

\newcommand{\hd}{\mathop{\rm pd}\nolimits}

\newcommand{\kernel}{\mathop{\rm Ker}\nolimits}

\newcommand{\assh}{\mathop{\rm Assh}\nolimits}
\newcommand{\minp}{\mathop{\rm Min}\nolimits}

\newcommand{\codim}{\mathop{\rm codim}\nolimits}

\renewcommand{\ker}{\mathop{\rm Ker}\nolimits}

\newcommand{\ol}{\overline}
\newcommand{\N}{\mathop{\rm N}}
\newcommand{\g}{\mathop{\rm G}\nolimits_0}
\newcommand{\k}[1]{\mathop{\rm K}\nolimits^{#1}_0}
\newcommand{\chow}[1]{\mathop{\rm A}\nolimits_{#1}}
\newcommand{\subq}{_{\Bbb Q}}
\newcommand{\chring}[1]{{\rm CH}^\cdot(#1)\subq}
\newcommand{\chnum}[1]{{\rm CH}^\cdot_{\rm num}(#1)\subq}
\newcommand{\chhom}[1]{{\rm CH}^\cdot_{\rm hom}(#1)\subq}

\newcommand{\p}{{\frak p}}
\newcommand{\q}{{\frak q}}
\newcommand{\m}{{\frak m}}
\newcommand{\n}{{\frak n}}

\begin{document}
\title{Numerical equivalence defined on Chow groups of Noetherian 
local rings}
\author{Kazuhiko Kurano (Meiji University)}
\date{}
\maketitle

\begin{abstract}
In the present paper, we define a notion of numerical equivalence 
on Chow groups or Grothendieck groups of Noetherian local rings,
which is an analogue of that on smooth projective varieties.
Under a mild condition, it is proved that
the Chow group modulo numerical equivalence is a finite dimensional 
${\Bbb Q}$-vector space, as in the case of smooth projective varieties.
Numerical equivalence on local rings is deeply related to 
that on smooth projective varieties.
For example, if Grothendieck's standard conjectures are true,
then a vanishing of Chow group (of local rings)
modulo numerical equivalence can be proven.

Using the theory of numerical equivalence,
the notion of numerically Roberts rings is defined.
It is proved that a Cohen-Macaulay local ring of positive characteristic
is a numerically Roberts ring if and only if
the Hilbert-Kunz multiplicity of a maximal primary ideal 
of finite projective dimension is always equal to its colength.
Numerically Roberts rings satisfy the vanishing property
of intersection multiplicities.

We shall prove another special case of the vanishing 
of intersection multiplicities using a vanishing of localized Chern
characters.
\end{abstract}

\section{Introduction}
For a smooth projective variety $X$ over a field,
the Chow ring $\chring{X}$ with rational coefficients is defined.
However, Chow rings are usually large and difficult to understand.
In order to study Chow rings,
we define numerical equivalence and sometimes consider
the Chow ring modulo numerical equivalence, denoted by $\chnum{X}$.
$\chnum{X}$ is a finite dimensional 
${\Bbb Q}$-vector space and a graded Gorenstein  ring.
Compared with the ordinary Chow ring $\chring{X}$,
$\chnum{X}$ is easier to study.

In the present paper, we define a notion of numerical equivalence 
on Chow groups or Grothendieck groups of Noetherian local rings 
(Definition~\ref{2.1}),
which is an analogue of that on smooth projective varieties.
We denote by $\ol{\chow{*}(R)\subq} = \oplus_i \ol{\chow{i}(R)\subq}$ 
(resp.\ $\ol{\g(R)\subq}$) 
the Chow group (resp.\ Grothendieck group) of a Noetherian local ring $R$
modulo numerical equivalence.
In many cases, the Chow group $\chow{*}(R)\subq$ 
and the Grothendieck group $\g(R)\subq$ are infinite dimensional 
${\Bbb Q}$-vector spaces.
The main theorem of the present paper is as follows:

\vspace{2mm}

\noindent
{\bf Theorem \ref{3.1}}
\newcommand{\maintheorem}{
Let $(R,\m)$ be a Noetherian excellent local ring that satisfies one of 
the following two conditions;
(1) $R$ contains ${\Bbb Q}$,
(2) $R$ is essentially of finite type over a field, ${\Bbb Z}$ or a complete 
discrete valuation ring.

Then, we have 
$\dim \ol{\k{\m}(R)\subq} = \dim \ol{\g(R)\subq} = \dim \ol{\chow{*}(R)\subq} < \infty$.
}
\begin{it}
\maintheorem
\end{it}

\vspace{2mm}

Numerical equivalence on Noetherian local rings is deeply related to that
on smooth projective varieties.  
(This will be discussed in Section~\ref{6}.)
In fact, let $A = \oplus_{n \geq 0}A_n$ be a standard graded ring over a field
$k = A_0$, and let $R = A_{A_+}$, where $A_+ = \oplus_{n > 0}A_n$.
Assume that $X = \proj(A)$ is smooth over $k$.
Let $h$ be the very ample divisor under the embedding.
It is known (cf.\ \cite{K11}) that there is an isomorphism
\[
\chring{X}/h \cdot \chring{X} \simeq \chow{*}(R)\subq .
\]
As will be shown in Section~\ref{6}, 
this isomorphism induces the natural surjection
\[
\chnum{X}/h \cdot \chnum{X} \longrightarrow \ol{\chow{*}(R)\subq} ,
\]
that is not always an isomorphism.
It is an isomorphism if and only if the natural map
\bangou
\begin{equation}\label{1.1}
\ker\left[\chring{X} \stackrel{h}{\rightarrow} \chring{X}\right] 
\rightarrow 
\ker\left[\chnum{X} \stackrel{h}{\rightarrow} \chnum{X}\right]
\end{equation}
is surjective, as will be shown in Section~\ref{6}.
(The map (\ref{1.1}) is studied in depth by Roberts and Srinivas~\cite{RS}.)
In particular, if the natural map $\chring{X} \longrightarrow \chnum{X}$
is an isomorphism, 
so is $\chow{*}(R)\subq \longrightarrow \ol{\chow{*}(R)\subq}$,
as shown in Theorem~\ref{coin}.
We shall see in Remark~\ref{6.9} that, 
if Grothendieck's standard conjectures are true, then 
\[
\mbox{
$\ol{\chow{j}(R)\subq} = 0$ for $j \leq \dim R/2$
} .
\]
This is equivalent to the condition that $\chi_{{\Bbb F}.}(M)$ is equal to 
$0$ for any ${\Bbb F}. \in C^\m(R)$ and any finitely generated $R$-module
$M$ with $\dim M \leq \dim R/2$.
We refer the reader to Section~\ref{2} for the definition of $C^\m(R)$ and 
$\chi_{{\Bbb F}.}(M)$.

In Section~\ref{sec3.5}, we shall study the invariant 
$\dim \ol{\g(R)\subq}$ of a local ring $R$.

Let $R$ be a homomorphic image of a regular local ring $S$.
Then, we have an isomorphism of ${\Bbb Q}$-vector spaces
\[
\tau_{R/S} : \g(R)\subq \longrightarrow \chow{*}(R)\subq
\]
by the singular Riemann-Roch theorem with base regular ring
$S$ (18.2 and 20.1 in \cite{F}).
Recall that the map $\tau_{R/S}$ is determined by both $R$ and $S$.
(The author does not know an explicit example such that 
$\tau_{R/S}$ actually depends on the choice of $S$.)
By Proposition~\ref{2.3}, we have an induced homomorphism
$\ol{\tau_{R/S}}$ that makes the following diagram commutative:
\[
\begin{array}{ccccc}
\g(R)\subq & \stackrel{\tau_{R/S}}{\longrightarrow} & \chow{*}(R)\subq & & \\
\downarrow 
& & 
 \downarrow & & 
\\
\ol{\g(R)\subq} & \stackrel{\ol{\tau_{R/S}}}{\longrightarrow} & 
\ol{\chow{*}(R)\subq} & = & \oplus_{i=0}^{\dim R}\ol{\chow{i}(R)\subq}
\end{array}
\]
We shall prove that the map $\ol{\tau_{R/S}}$ is independent 
of the  choice of $S$ in Section~\ref{4}.

Using the map $\ol{\tau_{R/S}}$,
the notion of numerically Roberts rings is defined in Section~\ref{5}.
Numerically Roberts rings are characterized in terms of  Dutta multiplicity or
Hilbert-Kunz multiplicity as follows:

\vspace{2mm}

\noindent
{\bf Theorem \ref{5.4}}
\newcommand{\DuttaHilbertKunz}{
Let $(R,\m)$ be a homomorphic image of a regular local ring.
\begin{itemize}
\item[(1)]
Then, $R$ is a numerically Roberts ring if and only if 
the Dutta multiplicity $\chi_\infty({\Bbb F}.)$ coincides with
the alternating sum of length of homology
\[
\chi({\Bbb F}.) = \sum_i(-1)^i\ell_R(H_i({\Bbb F}.))
\]
for any ${\Bbb F}. \in C^\m(R)$.
\item[(2)]
Assume that $R$ is a Cohen-Macaulay ring of characteristic $p$,
where $p$ is a prime number.
Then, $R$ is a numerically Roberts ring if and only if 
the Hilbert-Kunz multiplicity $e_{HK}(J)$ of $J$ coincides with the colength
$\ell_R(R/J)$ for any $\m$-primary ideal $J$ of finite projective dimension.
\end{itemize}
}
\begin{it}
\DuttaHilbertKunz
\end{it}

\vspace{2mm}

\noindent
Roberts rings defined in \cite{K16} are numerically Roberts rings.
As in Remark~\ref{newrem6.5}, numerically Roberts rings satisfy the
vanishing property of intersection multiplicities.
As in \cite{K16}, the category of Roberts rings contains
rings of dimension at most $1$, complete intersections,
quotient singularities, Galois extensions of regular local rings,
affine cones of abelian varieties, and many others.
As in Example~\ref{5.5}, the category of numerically Roberts rings contains
integral domains of dimension at most $2$,
Gorenstein rings of dimension $3$, and many others.

A vanishing of Chow groups modulo numerical equivalence will be discussed
in Section~\ref{7}.
We shall prove the following theorem:

\vspace{2mm}

\noindent
{\bf Theorem~\ref{7.1}}
\newcommand{\vanishing}{
Let $(R,\m)$ be a $d$-dimensional Noetherian local domain
that is a homomorphic image of a regular local ring.
Assume that there exists a regular alteration
$\pi : Z \rightarrow \spec(R)$.
Then, we have $\ol{\chow{t}(R)\subq} = 0$ for $t < d -\dim \pi^{-1}(\m)$.
}
\begin{it}
\vanishing
\end{it}

\vspace{2mm}

\noindent
By this theorem,  we have
\[
\dim \pi^{-1}(\m) \geq d - \min\{ t \mid \ol{\chow{t}(R)\subq} \neq 0 \} 
\]
for any resolution of singularities $\pi : Z \rightarrow \spec(R)$.

In Theorem~\ref{8.1}, we shall prove a special case of the vanishing property
of intersection multiplicities as follows:

\vspace{2mm}

\noindent
{\bf Theorem~\ref{8.1}}
\newcommand{\intersectionmultiplicity}{
Let $(R,\m)$ be a $d$-dimensional Noetherian local domain.
Assume that $R$ is a homomorphic image of an excellent regular local ring $S$
and there exists a regular alteration 
$\pi : Z \rightarrow \spec(R)$.
Let $Y$ be a closed subset of $\spec(R)$ such that
\[
\pi|_{Z \setminus \pi^{-1}(Y)} : Z \setminus \pi^{-1}(Y)
\rightarrow \spec(R) \setminus Y
\]
is finite.
If $\dim \pi^{-1}(Y) \leq d/2$,
then $R$ satisfies the vanishing property, that is,
\[
\sum_i(-1)^i\ell_R(\tor{i}{R}{M}{N}) = 0
\]
for finitely generated $R$-modules $M$ and $N$ 
such that 
$\hd_RM < \infty$, $\hd_RN < \infty$, 
$\ell_R(M\otimes_RN) < \infty$ and $\dim M + \dim N < d$,
where $\hd_R$ denotes the projective dimension as an $R$-module.
}
\begin{it}
\intersectionmultiplicity
\end{it}

\vspace{2mm}

In Theorem~\ref{8.5}, we give an another proof of the vanishing theorem
of the first localized Chern characters due to Roberts~\cite{Rmac}.

The author thanks the referee for many valuable comments.

\section{Numerical equivalence}\label{2}
Throughout the present paper, 
we assume that all local rings are homomorphic
images of regular local rings.
We denote by ${\Bbb Z}$, ${\Bbb Q}$ and ${\Bbb C}$ 
the ring of integers, the field of rational numbers and
the field of complex numbers, respectively.

Let $(R, \m)$ be a $d$-dimensional Noetherian 
local ring.
Let $C^\m(R)$ be the category of bounded complexes of finitely 
generated $R$-free modules and $R$-linear maps such that
its support is in $\{ \m \}$, i.e., each homology is of finite length
as an $R$-module.
A morphism in $C^\m(R)$ is a chain homomorphism of $R$-linear maps.
We sometimes call $C^\m(R)$ the category of bounded $R$-free 
complexes with support in $\{ \m \}$.

We define the {\em Grothendieck group} of complexes with support in $\{ \m \}$
as
\[
\k{\m}(R) = \left.
\bigoplus_{{\Bbb F}. \in C^\m(R)} {\Bbb Z} \cdot [{\Bbb F}.] \right/ P ,
\]
where 
$[{\Bbb F}.]$ is a free basis corresponding to the isomorphism class
containing a complex
\[
{\Bbb F}. \ : \ \cdots \rightarrow F_n \rightarrow F_{n-1} \rightarrow \cdots
\]
contained in $C^\m(R)$, and
$P$ is the subgroup generated by both
\[
\{ [{\Bbb G}.] - [{\Bbb F}.] - [{\Bbb H}.] \mid
\mbox{
$0 \rightarrow {\Bbb F}. \rightarrow {\Bbb G}. \rightarrow {\Bbb H}.
\rightarrow 0$ 
is exact in $C^\m(R)$} \}
\]
and 
\[
\{ [{\Bbb F}.] - [{\Bbb G}.] \mid 
\mbox{${\Bbb F}. \rightarrow {\Bbb G}.$ is a quasi-isomorphism
in $C^\m(R)$}
\} .
\]

We define the {\em Grothendieck group} of finitely generated $R$-modules
as
\[
\g(R) = 
\left.
\bigoplus_{M \in {\cal M}(R)} {\Bbb Z} \cdot [M] \right/ Q ,
\]
where ${\cal M}(R)$ is the category of finitely generated $R$-modules and 
$R$-linear maps, 
$[M]$ is a free basis corresponding to the isomorphism class
containing $M$, and
$Q$ is the subgroup generated by
\[
\{ [M] - [L] - [N] \mid
\mbox{
$0 \rightarrow L \rightarrow M \rightarrow N
\rightarrow 0$ 
is exact in ${\cal M}(R)$} \} .
\]

Let $\chow{*}(R) = \oplus_{i = 0}^d\chow{i}(R)$ be the {\em Chow group}
of the affine scheme $\spec(R)$, i.e.,
\[
\chow{i}(R) =  \left. \bigoplus {\Bbb Z} \cdot [\spec(R/\p)] \right/ 
{\rm Rat}_i(R) , 
\]
where the sum as above is taken over all prime ideals $\p$ of $\dim R/\p = i$,
$[\spec(R/\p)]$ is a free basis corresponding to a prime ideal $\p$,
and ${\rm Rat}_i(R)$ is the subgroup generated by rational equivalence 
(cf.\ Chapter~1 of \cite{F}).

For an additive group $A$, we denote by $A\subq$ the tensor product
$A \otimes_{\Bbb Z}{\Bbb Q}$.

Let $Z$ be a Noetherian scheme and let $X$ be a closed subset of $Z$.
We denote by $C^X(Z)$ the category of bounded complexes of vector bundles 
on $Z$ with support in $X$.
Let $\k{X}(Z)$ be 
the Grothendieck group of $C^X(Z)$.
Let $\g(Z)$ be the Grothendieck group of coherent sheaves on $Z$.
Let $\chow{*}(Z)$ be the Chow group of $Z$.
We refer the reader to Fulton~\cite{F}, Gillet-Soul\'e~\cite{GS} and
Srinivas~\cite{Sr} for definitions and basic properties on 
$\k{X}(Z)$, $\g(Z)$ and $\chow{*}(Z)$.
For ${\Bbb G}. \in C^X(Z)$, the {\em localized Chern characters}
(Chapter~18 of \cite{F}) 
\[
{\rm ch}({\Bbb G}.) = \sum_{i \geq 0}{\rm ch}_i({\Bbb G}.) :
\chow{*}(Z)\subq \longrightarrow \chow{*}(X)\subq
\]
are defined as operators on the Chow group, 
and for $\eta \in \chow{j}(Z)\subq$, we have
\[
{\rm ch}_i({\Bbb G}.)(\eta) \in \chow{j-i}(X)\subq
\]
for each $i$.

For ${\Bbb F}. \in C^\m(R)$ and $M \in {\cal M}(R)$,
we set
\[
\chi_{{\Bbb F}.}(M) = \sum_i(-1)^i\ell_R(H_i({\Bbb F}.\otimes_AM)) ,
\]
where $\ell_R( \ )$ denotes the {\em length} as an $R$-module.
We obtain well-defined maps
\bangou
\begin{equation}\label{eandv}
\begin{array}{l}
e : \k{\m}(R)\subq \otimes \g(R)\subq \longrightarrow {\Bbb Q} \\
v : \k{\m}(R)\subq \otimes \chow{*}(R)\subq \longrightarrow {\Bbb Q}
\end{array}
\end{equation}
that satisfy $e([{\Bbb F}.] \otimes [M] ) = \chi_{{\Bbb F}.}(M)$ and
$v([{\Bbb F}.] \otimes [\spec(R/\p)])
= {\rm ch}({\Bbb F}.)([\spec(R/\p)])$, where 
\[
{\rm ch}({\Bbb F}.) : \chow{*}(R)\subq \rightarrow \chow{*}(R/\m)\subq = 
{\Bbb Q} \cdot [\spec(R/\m)] = {\Bbb Q}
\]
is the localized Chern character of ${\Bbb F}. \in C^\m(R)$
(cf.\ Proposition~18.1 (b) and Example~18.1.4 in Fulton~\cite{F}).

We define cycles that are numerically equivalent to $0$ and
groups modulo numerical equivalence as follows.

\bangou
\begin{Definition}\label{2.1}
\begin{rm}
We define subgroups consisting of elements numerically equivalent to $0$
as 
\begin{eqnarray*}
\N\k{\m}(R)\subq & = & \{ \alpha \in \k{\m}(R)\subq \mid
\mbox{$e(\alpha \otimes \beta) = 0$ for any $\beta \in \g(R)\subq$} \} , \\
\N\g(R)\subq & = & \{ \beta \in \g(R)\subq \mid
\mbox{$e(\alpha \otimes \beta) = 0$ for any $\alpha \in \k{\m}(R)\subq$} \} , 
\\
\N\chow{i}(R)\subq & = & \{ \gamma \in \chow{i}(R)\subq \mid
\mbox{$v(\alpha \otimes \gamma) = 0$ for any $\alpha \in \k{\m}(R)\subq$} \} 
\end{eqnarray*}
for $i = 0, \ldots, d$.

We define groups modulo numerical equivalence as
\begin{eqnarray*}
\ol{\k{\m}(R)\subq} & = & \left. \k{\m}(R)\subq \right/ \N\k{\m}(R)\subq , \\
\ol{\g(R)\subq} & = & \left. \g(R)\subq \right/ \N\g(R)\subq , \\
\ol{\chow{i}(R)\subq} & = & \left. \chow{i}(R)\subq \right/ \N\chow{i}(R)\subq
\end{eqnarray*}
for $i = 0, \ldots, d$.
\end{rm}
\end{Definition}

By definition, $e$ induces a map
\[
\ol{e} : \ol{\k{\m}(R)\subq} \otimes \ol{\g(R)\subq} \longrightarrow {\Bbb Q} 
\]
that satisfies
$\ol{e}(\ol{\alpha} \otimes \ol{\beta}) = e(\alpha \otimes \beta)$ 
where $\ol{\alpha}$ and $\ol{\beta}$ denote
the images of $\alpha$ and $\beta$, respectively.

Let $R$ be a homomorphic image of a regular local ring $S$.
By the singular Riemann-Roch theorem with base regular ring $S$
(18.2 and 20.1 in \cite{F}),
we have an isomorphism of ${\Bbb Q}$-vector spaces
\[
\tau_{R/S} : \g(R)\subq \longrightarrow \chow{*}(R)\subq .
\]
Recall that the map $\tau_{R/S}$ as above is defined using
not only $R$ but also $S$.\footnote{
Let $Z$ be a regular scheme.
For a scheme $X$ of finite type over $Z$, we have an isomorphism
of ${\Bbb Q}$-vector spaces $\tau_{X/Z} : 
\g(X)\subq \longrightarrow \chow{*}(X)\subq$ by the singular 
Riemann-Roch theorem with regular base scheme $Z$ (20.1 in \cite{F}).

If $R$ is a homomorphic image of a regular local ring $S$,
we denote $\tau_{\spec(R)/\spec(S)}$ simply by $\tau_{R/S}$.

The map $\tau_{X/Z}$ usually depends on the choice of $Z$,
as will be shown in the following example.
Assume that $X$ is a smooth algebraic variety over a field $k$.
Set $Z = \spec(k)$.
Then, by 18.2 in \cite{F}, we have $\tau_{X/Z}([{\cal O}_X]) = 
{\rm td}(\Omega_{X/k}^\vee) = 1-c_1(K_X)/2 + \cdots$, that is, the Todd class 
(Example~3.2.4 in \cite{F}) of the tangent sheaf of $X$.
On the other hand, we have $\tau_{X/X}([{\cal O}_X]) = 1$ 
by definition.
Therefore, if $X$ is not ${\Bbb Q}$-Gorenstein,
then $\tau_{X/Z}([{\cal O}_X]) \neq \tau_{X/X}([{\cal O}_X])$.

However, in the case of $X = \spec(R)$ such that $R$ is a Noetherian local
ring, the author does not know any example such that $\tau_{X/Z}$ actually
depends on the choice of a regular base scheme $Z$ 
(cf.\ Section~4 in \cite{K16}).
}

Note that
\[
\N\k{\m}(R)\subq  =  \{ \alpha \in \k{\m}(R)\subq \mid
\mbox{$v(\alpha \otimes \gamma) = 0$ for any $\gamma \in \chow{*}(R)\subq$} \}
\]
because the diagram
\bangou
\begin{equation}\label{2.2}
\begin{array}{ccc}
\g(R)\subq & \stackrel{\tau_{R/S}}{\longrightarrow} & \chow{*}(R)\subq \\
{\scriptstyle \chi_{{\Bbb F}.}} \downarrow 
\hphantom{{\scriptstyle \chi_{{\Bbb F}.}}}
& & 
{\scriptstyle {\rm ch}({\Bbb F}.)} \downarrow 
\hphantom{{\scriptstyle {\rm ch}({\Bbb F}.)}} \\
{\Bbb Q} & = & {\Bbb Q} 
\end{array}
\end{equation}
is commutative for each ${\Bbb F}. \in C^\m(R)$ by the local Riemann-Roch 
theorem~(Example~18.3.12 in \cite{F}).

In order to construct an isomorphism $\ol{\tau_{R/S}}$ between 
$\ol{\g(R)\subq}$ and 
$\oplus_i\ol{\chow{i}(R)\subq}$,
we need the following proposition.

\bangou
\begin{Proposition}\label{2.3}
With notation as above, we have
\[
\tau_{R/S}(\N\g(R)\subq) = \bigoplus_{i = 0}^d \N\chow{i}(R)\subq .
\]
\end{Proposition}

\proof
Since diagram~(\ref{2.2})
is commutative for each ${\Bbb F}. \in C^\m(R)$, we have
\[
\tau_{R/S}(\N\g(R)\subq) = \{ \gamma \in \chow{*}(R)\subq \mid
\mbox{${\rm ch}({\Bbb F}.)(\gamma) = 0$ for any ${\Bbb F}. \in C^\m(R)$}
\} .
\]
By definition, it is easy to see that $\tau_{R/S}(\N\g(R)\subq) \supseteq 
\oplus_{i = 0}^d \N\chow{i}(R)\subq$.

We shall prove the opposite containment.
Let $\beta$ be an element of $\N\g(R)\subq$.
Set $\tau_{R/S}(\beta) = \gamma_d + \gamma_{d-1} + \cdots
+ \gamma_0$, where $\gamma_i \in \chow{i}(R)\subq$ for each $i$.
We want to show ${\rm ch}({\Bbb F}.)(\gamma_i) = 0$
for any ${\Bbb F}. \in C^\m(R)$ and any $i$.

For a positive integer $n$, we denote by 
$\psi^n : \k{\m}(R)\subq \rightarrow \k{\m}(R)\subq$ the 
$n$-th Adams operation defined by Gillet and Soul\'e~\cite{GS}.
For $n$, $i$ and ${\Bbb F}.$,
we have ${\rm ch}_i(\psi^n({\Bbb F}.)) = n^i{\rm ch}_i({\Bbb F}.)$
by Theorem~3.1 in \cite{K13}.
Then, we have
\begin{eqnarray*}
0 & = & \chi_{\psi^n({\Bbb F}.)}(\beta) \\
& = & {\rm ch}(\psi^n({\Bbb F}.))(\tau_{R/S}(\beta)) \\
& = & \sum_{i=0}^d{\rm ch}_i(\psi^n({\Bbb F}.))(\gamma_i) \\
& = & \sum_{i=0}^dn^{i} \cdot {\rm ch}_i({\Bbb F}.)(\gamma_i) \\
& = & \sum_{i=0}^dn^{i} \cdot {\rm ch}({\Bbb F}.)(\gamma_i) 
\end{eqnarray*}
for any positive integer $n$.
Therefore, we have ${\rm ch}({\Bbb F}.)(\gamma_i) = 0$ for each $i$.
\qed

We denote by $\ol{\chow{*}(R)\subq}$ the direct sum of 
$\ol{\chow{i}(R)\subq}$'s.
By the previous proposition,
we have an isomorphism $\ol{\tau_{R/S}}$ that makes the following diagram
commutative:
\[
\begin{array}{ccccc}
\g(R)\subq & \stackrel{\tau_{R/S}}{\longrightarrow} & \chow{*}(R)\subq & & \\
\downarrow 
& & 
 \downarrow & & 
\\
\ol{\g(R)\subq} & \stackrel{\ol{\tau_{R/S}}}{\longrightarrow} & 
\ol{\chow{*}(R)\subq} & = & \oplus_{i=0}^d\ol{\chow{i}(R)\subq}
\end{array}
\]
where the vertical maps are the natural projections.

The map $\tau_{R/S}$ is constructed using not only $R$ but also $S$.
However, it will be proved in Section~\ref{4} that
the map $\ol{\tau_{R/S}}$ is independent of the choice of $S$.

\bangou
\begin{Remark}\label{2.4}
\begin{rm}
For any local domain $(R,\m)$ of $\dim R \leq 2$, we shall show that
$\ol{\g(R)\subq}$ is spanned by $\ol{[R]}$ 
as a ${\Bbb Q}$-vector space by Proposition~\ref{3.5}.
If $R$ is an affine cone of a smooth curve over ${\Bbb C}$ of positive genus,
then $\dim \g(R)\subq = \infty$ as in Example~\ref{5.6}.
Therefore, $\N\g(R)\subq \neq 0$ in this case.

Let ${\Bbb K}.$ be the Koszul complex 
with respect to a system of parameters $\underline{a}$ of a local ring $R$.
Then, it is well known that $\chi_{{\Bbb K}.}([R])$ is 
the multiplicity of $R$ with respect to the ideal $(\underline{a})$.
Therefore, ${\Bbb K}.$ is not contained in $\N\k{\m}(R)\subq$,
and $\ol{\g(R)\subq}$ never coincides with $0$.
In particular, if $R$ is a regular local ring, then we have
$\g(R)\subq = \ol{\g(R)\subq} = {\Bbb Q}$.

For a finitely generated $R$-module $M$ with $\dim M < \dim R$, we have
\bangou
\begin{equation}\label{koszul}
\chi_{{\Bbb K}.}([M]) = 0 .
\end{equation}
However, for an arbitrary ${\Bbb F}. \in C^\m(R)$,
the equality (\ref{koszul}) does not always hold true
(cf.\ Dutta-Hochster-MacLaughlin~\cite{DHM}, Levin~\cite{L},
Miller-Singh~\cite{MS} and Roberts-Srinivas~\cite{RS}).
By Example~\ref{6.6}, 
we know the following:

Let $m$ and $n$ be positive integers such that $n \geq m \geq 2$.
Suppose
\[
R = \left. k[x_{ij} \mid \mbox{$i = 1, \ldots, m$; $j = 1, \ldots, n$}]_{
(x_{ij} \mid i, \ j)}
\right/
I_2(x_{ij}) ,
\]
where $I_2(x_{ij})$ is the ideal
generated by all the $2 \times 2$ minors of the $m \times n$ matrix $(x_{ij})$.
It is well known that the dimension of $R$ is $m + n -1$.
Then, for $s = n, n+1, \ldots, m + n -1$, 
there is a complex ${\Bbb H}\{ s \} . \in C^\m(R)$ 
that satisfies the following 
two properties:
\begin{enumerate}
\item
For any finitely generated $R$-module $M$ with $\dim M < s$,
$\chi_{{\Bbb H}\{ s \} .}([M]) = 0$.
\item
There exists a finitely generated $R$-module $N_s$ of dimension $s$ 
such that $\chi_{{\Bbb H}\{ s \} .}([N_s]) \neq 0$.
\end{enumerate}
\end{rm}
\end{Remark}

\section{Proof of the main theorem
}\label{3}
The aim of this section is to prove the following theorem:

\bangou
\begin{Theorem}\label{3.1}
\maintheorem
\end{Theorem}

Since $\ol{\tau_{R/S}}$ is an isomorphism,
we have $\dim \ol{\g(R)\subq} = \dim \ol{\chow{*}(R)\subq}$.
Furthermore, 
since the pairing
$\ol{e} : \ol{\k{\m}(R)\subq} \otimes \ol{\g(R)\subq} \longrightarrow {\Bbb Q}$
is perfect,
we have $\dim \ol{\k{\m}(R)\subq} = \dim \ol{\g(R)\subq}$
if $\dim \ol{\g(R)\subq} < \infty$.
Therefore, it is sufficient to prove $\dim \ol{\chow{*}(R)\subq} < \infty$.

\bangou
\begin{Lemma}\label{3.2}
Let $(A,\p) \stackrel{f}{\longrightarrow} (B,\q)$ be a finite morphism
of Noetherian local rings, that is, $B$ is a finitely generated
$A$-module.
We denote by $f_* : \chow{*}(B)\subq \rightarrow \chow{*}(A)\subq$ the 
induced map by the proper morphism $\spec(B) \rightarrow \spec(A)$.
Then, there is a map $\ol{f_*} : \ol{\chow{*}(B)\subq} \rightarrow
\ol{\chow{*}(A)\subq}$ that makes the following diagram commutative:
\[
\begin{array}{ccc}
\chow{*}(B)\subq & \stackrel{f_*}{\longrightarrow} & \chow{*}(A)\subq \\
\downarrow 
& & 
 \downarrow 
\\
\ol{\chow{*}(B)\subq} & \stackrel{\ol{f_*}}{\longrightarrow} & 
\ol{\chow{*}(A)\subq}
\end{array}
\]
where the vertical maps are the natural projections.
\end{Lemma}

\proof
It is sufficient to show $f_*(\alpha) \in \N\chow{i}(A)\subq$ for each $i$
and for each $\alpha \in \N\chow{i}(B)\subq$.

For ${\Bbb F}. \in C^\p(A)$, the complex ${\Bbb F}. \otimes_AB$ 
is contained in $C^\q(B)$ since the closed fibre of the morphism
$\spec(B) \rightarrow \spec(A)$ coincides with $\{ \q \}$.
By Definition~17.1 $({\rm C}_1)$ and Theorem~18.1 in \cite{F}, we have 
\begin{eqnarray*}
& & {\rm ch}({\Bbb F}.)(f_*(\alpha)) \\
& = & [B/\q : A/\p] \cdot {\rm ch}({\Bbb F}.\otimes B)(\alpha) \\
& = & 0 
\end{eqnarray*}
because $\alpha \in \N\chow{i}(B)\subq$.
Thus, we obtain $f_*(\alpha) \in \N\chow{i}(A)\subq$.
\qed

Let $(R,\m)$ be a local ring that satisfies 
the assumptions in Theorem~\ref{3.1}.
Let $\{ \p_1, \ldots, \p_r \}$ be the set of minimal prime ideals 
of $R$.
Let $R_i = R/\p_i$, and let $f_i : R \rightarrow R_i$ denote the
 projection for $i = 1, \ldots, r$.
By the previous proposition,
we have the following commutative diagram:
\[
\begin{array}{ccc}
\oplus_i\chow{*}(R_i)\subq & \stackrel{\sum_i{f_i}_*}{\ya{aaaaaaaa}} & 
\chow{*}(R)\subq \\
\downarrow 
& & 
 \downarrow 
\\
\oplus_i\ol{\chow{*}(R_i)\subq} & \stackrel{\sum_i\ol{{f_i}_*}}{\ya{aaaaaaaa}} 
& 
\ol{\chow{*}(R)\subq}
\end{array}
\]
where all maps are surjections.
Therefore, we have only to show that $\dim \ol{\chow{*}(R_i)\subq} < \infty$
for each $i$.

Hence, we may assume that the given local ring $(R,\m)$ is an integral domain.

Then, by Hironaka~\cite{H} or de~Jong~\cite{dJ},
there exists a projective regular alteration $\pi : Z \rightarrow \spec(R)$,
that is, a projective generically finite morphism such that $Z$ is 
a regular scheme.
Furthermore, 
we may assume that $\pi^{-1}(\m)_{\rm red}$ is a simple normal crossing
divisor.
Let $\pi^{-1}(\m)_{\rm red} = E_1 \cup \cdots \cup E_t$ be the 
irreducible decomposition.
Then, each $E_l$ is a regular projective variety over $R/\m$
with ${\rm codim}_ZE_l = 1$ for each $l = 1, \ldots, t$.

We denote by $\chnum{E_l}$ the Chow group of $E_l$ (with rational coefficients)
modulo numerical equivalence (Chapter~19 in \cite{F}).
Then, $\dim \chnum{E_l}$ is finite as in Example~19.1.4 in \cite{F}.
(In Example~19.1.4,  it is assumed that the base field is algebraically closed.
However, it is easy to remove this assumption.)
Then, Theorem~\ref{3.1} follows from the following claim:

\bangou
\begin{Claim}\label{3.3}
With notation as above,
$\ol{\chow{*}(R)\subq}$ is a subquotient of
$\oplus_{l=1}^t\chnum{E_l}$.
\end{Claim}

The claim is proven as follows.

Since $\pi_* : \chow{*}(Z)\subq \rightarrow \chow{*}(R)\subq$ is surjective,
we have a map $s : \chow{*}(R)\subq \rightarrow \chow{*}(Z)\subq$
of ${\Bbb Q}$-vector spaces such that
$\pi_* \cdot s = 1$.

For each $l = 1, \ldots, t$,
we denote by $j_l$ the inclusion $E_l \rightarrow Z$.
We denote by $\varphi$ the composite map 
\[
\chow{*}(R)\subq \stackrel{s}{\longrightarrow}
\chow{*}(Z)\subq \stackrel{\sum j_l^*}{\ya{aaaaaa}}
\oplus_l\chow{*}(E_l)\subq \longrightarrow
\oplus_l\chnum{E_l} ,
\]
where the last map is the natural projection and
$j_l^*$ is the map 
taking the intersection with the effective Cartier divisor $E_l$
(see Chapter~2 in \cite{F}).

In order to prove the claim,
it is sufficient to show that the kernel of $\varphi$ is contained
in $\oplus_i\N\chow{i}(R)\subq$.

Assume that $\gamma$ is an element of $\chow{*}(R)\subq$ 
such that $\varphi(\gamma) = 0$.
We shall prove ${\rm ch}({\Bbb F}.)(\gamma) = 0$ for any
${\Bbb F}. \in C^\m(R)$.
Note that, since $\varphi(\gamma) = 0$,
\bangou
\begin{equation}\label{abc}
\mbox{$j_l^* \cdot s(\gamma)$ is equal to $0$ in $\chnum{E_l}$ for each $l$.}
\end{equation}

Set $E = \pi^{-1}(\m)$.
Since the diagram 
\[
\begin{array}{ccc}
E & \longrightarrow & Z \\
{\scriptstyle \pi'} \downarrow \phantom{{\scriptstyle \pi'}}
& & 
{\scriptstyle \pi} \downarrow \phantom{{\scriptstyle \pi}} \\
\spec(R/\m) & \rightarrow & \spec(R)
\end{array}
\]
is a fibre square, we have
\begin{eqnarray*}
{\rm ch}({\Bbb F}.)(\gamma) & = & 
{\rm ch}({\Bbb F}.)(\pi_* \cdot s (\gamma)) \\
& = & \pi'_* {\rm ch}(\pi^*{\Bbb F}.)(s (\gamma)) .
\end{eqnarray*}
Note that $\pi^*{\Bbb F}.$ is a complex in $C^E(Z)$,
that is, $\pi^*{\Bbb F}.$ is a bounded complex of vector bundles on $Z$
with support in $E$.
Consider the following commutative diagram:
\[
\begin{array}{ccc}
\oplus_l\k{E_l}(Z)\subq & \stackrel{\sum{i_l}_*}{\ya{aaaaaaaa}} & \k{E}(Z)\subq \\
{\scriptstyle \oplus\chi} \downarrow \phantom{\scriptstyle \oplus\chi} & & 
{\scriptstyle \chi} \downarrow \phantom{\scriptstyle \chi} \\
\oplus_l\g(E_l)\subq & \stackrel{\sum{i_l}_*}{\ya{aaaaaaaa}} & \g(E)\subq 
\end{array}
\]
where $\chi : \k{Y}(Z)\subq \rightarrow \g(Y)\subq$ is a map defined by
$\chi({\Bbb G}.) = \sum_i(-1)^i[H_i({\Bbb G}.)]$
for each closed subset $Y$ of $Z$ and for each ${\Bbb G}. \in C^Y(Z)$, 
${i_l}_* : \k{E_l}(Z)\subq \rightarrow \k{E}(Z)\subq$ is a map given by
${i_l}_*([{\Bbb H}.]) = [{\Bbb H}.]$ for ${\Bbb H}. \in C^{E_l}(Z)$, and
${i_l}_* : \g(E_l)\subq \rightarrow \g(E)\subq$ is a map given by 
${i_l}_*([{\cal F}]) = [{\cal F}]$ for each coherent ${\cal O}_{E_l}$-module
${\cal F}$.
Since $Z$ is a regular scheme, the vertical maps are isomorphisms
by Lemma~1.9 in \cite{GS}.
The bottom map is surjective,
because $E_{\rm red} = E_1 \cup \cdots \cup E_t$.
Therefore, the top map is also surjective
and there exist $\delta_l \in \k{E_l}(Z)\subq$ for $l = 1, \ldots, t$
such that
\[
[\pi^*{\Bbb F}.] = \sum_l {i_l}_*(\delta_l) \ \ \mbox{in $\k{E}(Z)\subq$} .
\]
For each $l$, $\pi_l : E_l \rightarrow \spec(R/\m)$ denotes the 
structure morphism.
Then, we have
\[
\pi'_* {\rm ch}(\pi^*{\Bbb F}.)(s (\gamma))
= \sum_l {\pi'}_*{\rm ch}({i_l}_*(\delta_l))(s(\gamma))
= \sum_l {\pi'}_*{i_l}_*{\rm ch}(\delta_l)(s(\gamma))
= \sum_l {\pi_l}_*{\rm ch}(\delta_l)(s(\gamma)) 
\]
since ${\rm ch}({i_l}_*(\delta_l)) = {i_l}_*{\rm ch}(\delta_l)$.

Here, we shall prove ${\pi_l}_*{\rm ch}(\delta_l)(s(\gamma)) = 0$
for each $l$.
We denote by $g_l$ the composite map 
\[
\k{E_l}(Z)\subq \stackrel{\sim}{\rightarrow}
 \g(E_l)\subq \stackrel{\sim}{\rightarrow} \k{E_l}(E_l)\subq .
\]
Set $\epsilon_l = g_l(\delta_l)$ for $l = 1, \ldots, t$.
Then, by Corollary~18.1.2 in \cite{F}, we have
\[
{\pi_l}_*{\rm ch}(\delta_l)(s(\gamma))
= {\pi_l}_*\left(
{\rm ch}(\epsilon_l) \cdot 
{\rm td}(j_l^*{\cal O}(E_l))^{-1} \cdot j_l^*(s(\gamma))
\right) =0 
\]
since $j_l^* \cdot s(\gamma)$ is equal to $0$ in $\chnum{E_l}$ as (\ref{abc}).

We have obtained ${\rm ch}({\Bbb F}.)(\gamma) = 0$.

We have completed the proof of Theorem~\ref{3.1}.

\bangou
\begin{Remark}\label{3.35}
\begin{rm}
Let $(R,\m)$ be a $d$-dimensional local ring that satisfies the assumptions
as in Theorem~\ref{3.1}.
Set
\begin{eqnarray*}
\N\k{\m}(R) & = & \{ \alpha \in \k{\m}(R) \mid
\mbox{$e(\alpha \otimes \beta) = 0$ for any $\beta \in \g(R)$} \} , \\
\N\g(R) & = & \{ \beta \in \g(R) \mid
\mbox{$e(\alpha \otimes \beta) = 0$ for any $\alpha \in \k{\m}(R)$} \} , 
\\
\N\chow{*}(R) & = & \{ \gamma \in \chow{*}(R) \mid
\mbox{$v(\alpha \otimes \gamma) = 0$ for any $\alpha \in \k{\m}(R)$} \} ,
\end{eqnarray*}
where $e$ and $v$ are the maps as in (\ref{eandv}).
It is easy to see that $\left( \k{\m}(R)/\N\k{\m}(R) \right)\subq$, 
$\left( \g(R)/\N\g(R) \right)\subq$ and 
$\left( \chow{*}(R)/\N\chow{*}(R) \right)\subq$ coincide with
$\ol{\k{\m}(R)\subq}$, $\ol{\g(R)\subq}$ and 
$\ol{\chow{i}(R)\subq}$, respectively.
Therefore, by Theorem~\ref{3.1}, $\k{\m}(R)/\N\k{\m}(R)$, 
$\g(R)/\N\g(R)$ and $\chow{*}(R)/\N\chow{*}(R)$ are torsion-free abelian 
groups of finite rank.
Note that the pairing
\[
\ol{e} : \left( \k{\m}(R)/\N\k{\m}(R) \right)\subq \otimes\subq
\left( \g(R)/\N\g(R) \right)\subq \longrightarrow {\Bbb Q}
\]
is perfect and it satisfies that $\ol{e}(a \otimes b) \in {\Bbb Z}$
for any $a \in \k{\m}(R)/\N\k{\m}(R)$ and any $b \in \g(R)/\N\g(R)$.

Suppose that $L_1$ and $L_2$ are torsion-free abelian groups of finite rank
with a perfect pairing $p : {L_1}\subq \otimes {L_2}\subq \rightarrow {\Bbb Q}$
such that $p(a \otimes b) \in {\Bbb Z}$ 
for any $a \in L_1$ and any $b \in L_2$.
Then, $L_1$ and $L_2$ can easily be proven to be finitely generated 
free abelian groups.

Based on this fact, we know that $\k{\m}(R)/\N\k{\m}(R)$ and 
$\g(R)/\N\g(R)$ are finitely generated free abelian groups.
By 18.1 (14) in Fulton~\cite{F},
we have ${\rm ch}({\Bbb F}.)(\gamma) \in {\Bbb Z}/d!$ for any
${\Bbb F}. \in C^{\m}(R)$ and $\gamma \in \chow{*}(R)$.
Thus, we can prove similarly that $\chow{*}(R)/\N\chow{*}(R)$ is also
a finitely generated free abelian group.
\end{rm}
\end{Remark}

\bangou
\begin{Remark}\label{3.4}
\begin{rm}
Let $(R,\m)$ be a $d$-dimensional Noetherian local ring.
By definition, for ${\Bbb F}. \in C^\m(R)$ and for
a non-negative integer $n$,
the following two conditions are equivalent:
\begin{enumerate}
\item
${\rm ch}_i({\Bbb F}.) : \chow{i}(R)\subq \rightarrow {\Bbb Q}$
is equal to $0$ for $i = 0, \ldots, n$, 
\item
$\chi_{{\Bbb F}.}(M)$ is equal to $0$
for any finitely generated 
$R$-module $M$ with $\dim M \leq n$.
\end{enumerate}

A Koszul complex of a system of parameters satisfies
the conditions as above with $n = d-1$.

%
%
On the other hand, 
for $s = n, n+1, \ldots, m + n -1$,
the complex ${\Bbb H}\{ s \}.$ in Remark~\ref{2.4}
satisfies following two conditions:
(1) ${\rm ch}_i({\Bbb H}\{ s \}.) : \chow{i}(R)\subq \rightarrow {\Bbb Q}$
is equal to $0$ for $i < s$, and
(2) ${\rm ch}_s({\Bbb H}\{ s \}.) \neq 0$.
\end{rm}
\end{Remark}

\bangou
\begin{Proposition}\label{3.5}
Let $(R,\m)$ be a $d$-dimensional Noetherian local ring.
\begin{itemize}
\item[(1)]
If $d > 0$, then we have $\chow{0}(R)\subq = \ol{\chow{0}(R)\subq} = 0$.
\item[(2)]
Assume that $\dim R/\p$ is at least $2$ for each minimal prime ideal 
$\p$ of $R$.
Then, $\ol{\chow{1}(R)\subq} = 0$.
\item[(3)]
The natural map  $\chow{d}(R)\subq \rightarrow \ol{\chow{d}(R)\subq}$
is an isomorphism.
\end{itemize}
\end{Proposition}

\proof
By definition, if $d > 0$, then we have $\chow{0}(R)\subq = 0$.
Since $\ol{\chow{0}(R)\subq}$ is a homomorphic image of $\chow{0}(R)\subq$,
we have $\ol{\chow{0}(R)\subq} = 0$ in this case.

Roberts~\cite{Rmac} proved ${\rm ch}_1({\Bbb F}.) = 0$
for any ${\Bbb F}. \in C^\m(R)$ under the assumption in (2) as above.
Therefore, we have $\ol{\chow{1}(R)\subq} = 0$ in this case.
(We shall give an another proof of the vanishing theorem of
the first localized Chern characters due to Roberts in 
Theorem~\ref{8.5}.)

We shall prove (3).
Set $\assh R = \{ \p_1, \ldots, \p_t \}$, that is
the set of minimal prime ideals of coheight $d$.
Then, $\chow{d}(R)\subq$ is the ${\Bbb Q}$-vector space 
with basis $\{ [\spec(R/\p_j)] \mid j = 1, \ldots, t \}$.
Suppose that $\sum_jn_j[\spec(R/\p_j)] \in \N\chow{d}(R)\subq$.
We want to show that $n_1 = \cdots = n_t = 0$.

For each $j = 1, \ldots, t$, take 
\[
x_j \in \left( \bigcap_{\uetuki{\q \in \minp(R)}{\q\neq \p_j}} \q  \right) 
\setminus \p_j .
\]
Then, $\sum_kx_k \not\in \p_j$ for $j = 1, \ldots, t$.
Therefore, we can choose $y_2, \ldots, y_d \in \m$ such that
$\sum_kx_k, y_2, \ldots, y_d$ is a system of parameters for $R$.
Then, it is easy to check that
$\sum_kx_k^{s_k}, y_2, \ldots, y_d$ is also a system of parameters for $R$
for any positive integers $s_1, \ldots, s_t$.
Let ${\Bbb K}(\sum_kx_k^{s_k}, y_2, \ldots, y_d).$ be the Koszul complex with respect
to $\sum_kx_k^{s_k}, y_2, \ldots, y_d$.

Since $R$ is a homomorphic image of a regular local ring $S$,
we have an isomorphism of ${\Bbb Q}$-vector spaces
\[
\tau_{R/S} : \g(R)\subq \longrightarrow \chow{*}(R)\subq 
\]
by the singular Riemann-Roch theorem~\cite{F}.
By the top term property~(Theorem~18.3~(5) in \cite{F}), we have
\[
\tau_{R/S}([R/\p_j]) = [\spec(R/\p_j)] + \gamma_{j, d-1} + \cdots + 
\gamma_{j 0},
\]
where $\gamma_{ji} \in \chow{i}(R)\subq$.
Then, by the local Riemann-Roch formula~(Example~18.3.12 in \cite{F}),
we have
\begin{eqnarray*}
\chi_{{\Bbb K}(\sum_kx_k^{s_k}, y_2, \ldots, y_d).}([R/\p_j])
& = & {\rm ch}_d({\Bbb K}(\sum_kx_k^{s_k}, y_2, \ldots, y_d).)([\spec(R/\p_j)])
\\ & & 
+ \sum_{i = 0}^{d-1}{\rm ch}_i({\Bbb K}(\sum_kx_k^{s_k}, y_2, \ldots, y_d).)
(\gamma_{ji}) .
\end{eqnarray*}
On the other hand, 
we have 
\[
{\rm ch}_i({\Bbb K}(\sum_kx_k^{s_k}, y_2, \ldots, y_d).)(\gamma_{ji}) = 0
\]
for $i = 0, \ldots, d-1$ by Remark~\ref{3.4}.
Therefore, we have
\[
\chi_{{\Bbb K}(\sum_kx_k^{s_k}, y_2, \ldots, y_d).}(\sum_jn_j[R/\p_j])
= {\rm ch}_d({\Bbb K}(\sum_kx_k^{s_k}, y_2, \ldots, y_d).)
(\sum_jn_j[\spec(R/\p_j)]) = 0
\]
since $\sum_jn_j[\spec(R/\p_j)] \in \N\chow{d}(R)\subq$.

By Lech's Lemma and a theorem of Auslander-Buchsbaum
(cf.\ p110, p111 in Matsumura~\cite{Ma}), we have
\begin{eqnarray*}
0 & = & 
\chi_{{\Bbb K}(\sum_kx_k^{s_k}, y_2, \ldots, y_d).}(\sum_jn_j[R/\p_j]) \\
& = & \sum_jn_j \chi_{{\Bbb K}(\sum_kx_k^{s_k}, y_2, \ldots, y_d).}([R/\p_j]) \\
& = & \sum_jn_j \chi_{{\Bbb K}(x_j^{s_j}, y_2, \ldots, y_d).}([R/\p_j]) \\
& = & \sum_jn_j e((x_j^{s_j}, y_2, \ldots, y_d), R/\p_j) \\
& = & \sum_jn_j s_j \cdot e((x_j, y_2, \ldots, y_d), R/\p_j) 
\end{eqnarray*}
for any $s_1, \ldots, s_t > 0$,
where $e(\ , \ )$ denotes the multiplicity.
Since $e((x_j, y_2, \ldots, y_d), R/\p_j) > 0$ for each $j$,
we have $n_1 = \cdots = n_t = 0$.
\qed

\section{Dimension of $\ol{\g(R)\subq}$ as a ${\Bbb Q}$-vector space}
\label{sec3.5}

By Theorem~\ref{3.1}, we know that the dimension of $\ol{\g(R)\subq}$ 
as a ${\Bbb Q}$-vector space is finite for a local ring $R$ that satisfies
a mild condition.
It is expected that $\dim \ol{\g(R)\subq} < \infty$ holds
for an arbitrary local ring $R$.

In this section, we study the dimension of $\ol{\g(R)\subq}$ as a 
${\Bbb Q}$-vector space.

As was shown in Remark~\ref{2.4},
$\dim \ol{\g(R)\subq}$ is always positive.
More precisely, by Theorem~\ref{3.1} and Proposition~\ref{3.5}~(3),
$\dim \ol{\g(R)\subq}$ is at least the number of prime ideals
contained in $\assh R$.
If $R$ is a regular local ring, then $\dim \ol{\g(R)\subq}$ is equal to one
since $\g(R)\subq = {\Bbb Q}$.

For a local domain $R$ of dimension $\leq 2$,
$\dim \ol{\g(R)\subq}$ is equal to $1$ by Proposition~\ref{3.5}.

First of all, we give some examples.

\bangou
\begin{Example}\label{3.5.1}
\begin{rm}
Let $m$ and $n$ be positive integers such that $n \geq m \geq 2$.
Suppose
\[
R = \left. k[x_{ij} \mid \mbox{$i = 1, \ldots, m$; $j = 1, \ldots, n$}]_{
(x_{ij} \mid i, \ j)}
\right/
I_2(x_{ij}) 
\]
as in Remark~\ref{2.4}.
In Example~\ref{6.6}, we will show that $\dim \ol{\g(R)\subq}$ 
is equal to $m$.

Let $X$ be the blowing-up of the projective space
${\Bbb P}^n_{\Bbb C}$ at $r$ distinct points,
where we suppose that $n \geq 2$.
Let $R$ be the local ring (at the homogeneous maximal ideal) of an affine
cone of $X$.
Then, by Example~\ref{6.5.5}, $\dim \ol{\g(R)\subq}$ is equal to $r+1$.

Hence, $\dim \ol{\g(R)\subq}$ has no upper bound even if $\dim R = 3$.
\end{rm}
\end{Example}

In the remainder of this section, we compare $\dim \ol{\g(A)\subq}$ with 
$\dim \ol{\g(B)\subq}$ for a homomorphism $A \rightarrow B$
of Noetherian local rings.

\bangou
\begin{Lemma}\label{4.1}
Let $f : (A,\p) \rightarrow (B,\q)$ be a flat local homomorphism of 
Noetherian local rings such that $B/\q$ is a finite algebraic extension
of $A/\p$.
\begin{itemize}
\item[(1)]
The map $f^* : \g(A)\subq \rightarrow \g(B)\subq$ induces
a map $\ol{f^*} : \ol{\g(A)\subq} \rightarrow \ol{\g(B)\subq}$ that makes
the following diagram commutative:
\[
\begin{array}{ccc}
\g(A)\subq & \stackrel{f^*}{\ya{aaaa}} & \g(B)\subq \\
\downarrow & & \downarrow \\
\ol{\g(A)\subq} & \stackrel{\ol{f^*}}{\ya{aaaa}} & \ol{\g(B)\subq}
\end{array}
\]
where the vertical maps are the projections.
\item[(2)]
Assume $\sqrt{\p B} = \q$.
Then, the map $\ol{f^*} : \ol{\g(A)\subq} \rightarrow \ol{\g(B)\subq}$ is 
injective.

In particular, $\dim \ol{\g(A)\subq} \leq \dim \ol{\g(B)\subq}$
in this case.
\end{itemize}
\end{Lemma}

\proof
Recall that $f^* : \g(A) \rightarrow \g(B)$ is a map defined by 
$f^*([M]) = [M \otimes_AB]$ for each finitely generated $A$-module $M$.

First, we prove (1).
It is sufficient to show that, for any $c \in \N\g(A)\subq$, 
$f^*(c)$ is contained in $\N\g(B)\subq$.

It is sufficient to show that, for any ${\Bbb F}. \in C^\q(B)$,
$\chi_{{\Bbb F}.}(f^*(c)) = 0$.

First we shall prove that there exist ${\Bbb G}. \in C^\p(A)$ and a chain map 
$\varphi : {\Bbb G}. \rightarrow {\Bbb F}.$ of $A$-linear maps
such that $\varphi$ is a quasi-isomorphism.
Recall that each homology of ${\Bbb F}. \in C^\q(B)$ is of finite length as an
$A$-module since $B/\q$ is a finite algebraic extension of $A/\p$.
By killing homology modules of ${\Bbb F}.$,
there exists an $A$-free complex ${\Bbb G}.$ and a chain map 
$\phi : {\Bbb G}. \rightarrow {\Bbb F}.$ of $A$-linear maps such that
\begin{enumerate}
\item
each $G_i$ is a finitely generated free $A$-module,
\item
${\Bbb G}.$ is bounded below, i.e., $G_i = 0$ for $i << 0$, and
\item
$\phi : {\Bbb G}. \rightarrow {\Bbb F}.$ is a quasi-isomorphism.
\end{enumerate}
Furthermore, we may assume that ${\Bbb G}.$ is a minimal complex,
that is, all boundaries of ${\Bbb G}.\otimes_AA/\p$ are $0$.
Note that each homology of ${\Bbb G}.$ is of finite length 
as an $A$-module.
We want to show that ${\Bbb G}.$ is bounded.
Since $\phi : {\Bbb G}. \rightarrow {\Bbb F}.$ is a quasi-isomorphism
of complexes of flat $A$-modules which are bounded below,
$\phi\otimes 1 : {\Bbb G}.\otimes_AA/\p \rightarrow 
{\Bbb F}.\otimes_AA/\p$ is also a quasi-isomorphism.
Therefore, we have
\[
(A/\p)^{\rank_AG_i} \simeq H_i({\Bbb G}.\otimes_AA/\p) \simeq
H_i({\Bbb F}.\otimes_AA/\p)
\]
for any $i$.
Since ${\Bbb F}.$ is bounded, so is ${\Bbb G}.$.

Since $\phi : {\Bbb G}. \rightarrow {\Bbb F}.$ is a quasi-isomorphism
of bounded complexes of flat $A$-modules,
$\phi\otimes 1 : {\Bbb G}.\otimes_AM \rightarrow {\Bbb F}.\otimes_AM$
is a quasi-isomorphism for any $A$-module $M$.
Since ${\Bbb F}.\otimes_AM = {\Bbb F}.\otimes_B(B \otimes_AM)$,
$H_i({\Bbb G}.\otimes_AM)$ is isomorphic to 
$H_i({\Bbb F}.\otimes_B(B \otimes_AM))$ as an $A$-module.
Therefore, we have
\[
\chi_{{\Bbb F}.}(f^*(c)) = 
\frac{1}{[B/\q : A/\p]}\chi_{{\Bbb G}.}(c) = 0
\]
since $c \in \N\g(A)\subq$.
We have completed the proof of (1) in Lemma~\ref{4.1}.

Next, we shall prove (2).
Assume that $\beta \in \g(A)\subq$ satisfies $f^*(\beta) \in \N\g(B)\subq$.
We want to show that $\beta \in \N\g(A)\subq$.
Let ${\Bbb G}.$ be a complex contained in $C^\p(A)$.
Since $\sqrt{\p B} = \q$, ${\Bbb G}.\otimes_AB$ is contained in $C^\q(B)$.
Then, we have
\[
\chi_{{\Bbb G}.}(\beta) = 
\frac{1}{\ell_B(B/\p B)}\chi_{{\Bbb G}.\otimes_AB}(f^*(\beta)) = 0 .
\]
Therefore, we have $\beta \in \N\g(A)\subq$.
\qed

\bangou
\begin{Theorem}\label{3.5.2}
Let $f : (A,\p) \rightarrow (B,\q)$ be a local homomorphism 
of Noetherian local rings.
\begin{itemize}
\item[(1)]
If $f$ is finite and injective, 
then $\dim \ol{\g(A)\subq} \leq \dim \ol{\g(B)\subq}$.
\item[(2)]
Assume that $A$ is an excellent normal local ring.
Suppose that the Frobenius map $A \rightarrow A$ is finite
if $A$ is of positive characteristic.
Let $K$ be a finite normal (algebraic) extension of $Q(A)$,
where $Q( \ )$ is the field of fractions.
Assume that $B$ is a local ring at a maximal ideal 
of the integral closure of $A$
in $K$.
Then, we have $\dim \ol{\g(A)\subq} \leq \dim \ol{\g(B)\subq}$.
\end{itemize}
\end{Theorem}

\proof
Assume that $f$ is finite and injective.
Then, by the lying-over theorem,
$f_* : \chow{*}(B)\subq \rightarrow \chow{*}(A)\subq$ is surjective.
By Lemma~\ref{3.2}, we obtain a surjective map
$\ol{f_*} : \ol{\chow{*}(B)\subq} \rightarrow \ol{\chow{*}(A)\subq}$. 
Therefore, by Theorem~\ref{3.1}, we have 
$\dim \ol{\g(A)\subq} \leq \dim \ol{\g(B)\subq}$.

Next, we shall prove (2).
Let $L$ be the intermediate field such that $Q(B)/L$ is a Galois extension
and $L/Q(A)$ is a purely inseparable extension.
Let $C$ (resp.\ $D$) be the integral closure of $A$ in $L$ (resp.\ $Q(B)$).
Since $A$ is excellent, both $C$ and $D$ are finite $A$-modules.

Suppose that $L \neq Q(A)$.
Then, the characteristic of $A$ is a prime number $p$ and, 
by the assumption, the Frobenius map $A \rightarrow A$ is finite.
Therefore, for a sufficiently large $e$, 
\[
A \subset C \subset A^{1/p^e} \simeq A .
\]
Since the maps as above are finite, we have
\[
\dim \ol{\g(A)\subq} \leq \dim \ol{\g(C)\subq} \leq 
\dim \ol{\g(A^{1/p^e})\subq} = \dim \ol{\g(A)\subq} 
\]
by (1).
Therefore, $\dim \ol{\g(A)\subq}$ coincides with $\dim \ol{\g(C)\subq}$.
Replacing $C$ with $A$, we may assume that $Q(B)$ is a Galois extension
of $Q(A)$.

Let $\m_1, \ldots, \m_t$ be the set of maximal ideals of $D$.
Suppose $B = D_{\m_1}$.
Let $G$ be the Galois group of the Galois extension $Q(B)/Q(A)$.
Let $H$ be the splitting group of $\m_1$, that is
\[
H = \{ \sigma \in G \mid \sigma(\m_1) = \m_1 \} .
\]
Set $K = Q(B)^H$.
Let $E$ be the integral closure of $A$ in $K$.
Set $\n = \m_1 \cap E$.
Since the induced map $E_{\n} \rightarrow D_{\m_1}$ is finite injective
(e.g., (41.2)~(1) in Nagata~\cite{LR}),
we have
\[
\dim \ol{\g(E_{\n})\subq} \leq \dim \ol{\g(D_{\m_1})\subq} .
\]
By (41.2) (2), (3) in \cite{LR}, we obtain
$\p E_{\n} = \n E_{\n}$ and $A/\p = E_{\n}/\n E_{\n}$.
Then, by (43.1) in \cite{LR}, 
$E_{\n}$ is flat over $A$.
Then, by Lemma~\ref{4.1}, we have
\[
\dim \ol{\g(A)\subq} \leq \dim \ol{\g(E_{\n})\subq} .
\]
The assertion (2) follows immediately from the above inequalities.
\qed

\section{$\ol{\tau_{R/S}}$ is independent of $S$}\label{4}

This section is devoted to proving that the map
$\ol{\tau_{R/S}} : \ol{\g(R)\subq} \rightarrow \ol{\chow{*}(R)\subq}$ is independent
of the choice of $S$.

The essential point is to prove that
the induced map $\ol{\g(R)\subq} \rightarrow \ol{\g(\hat{R})\subq}$ by 
completion $R \rightarrow \hat{R}$ is injective as follows.

\bangou
\begin{Theorem}\label{4.2}
\begin{itemize}
\item[(1)]
Let $R$ be a Noetherian local ring.
Then, the induced map $\ol{\g(R)\subq} \rightarrow \ol{\g(\hat{R})\subq}$ 
by the completion $R \rightarrow \hat{R}$ is injective.
\item[(2)]
Let $R$ be a homomorphic image of a regular local ring $S$.
Then, the induced map $\ol{\tau_{R/S}} : \ol{\g(R)\subq} \rightarrow
\ol{\chow{*}(R)\subq}$ is independent of the choice of $S$.
\end{itemize}
\end{Theorem}

\proof
The first assertion follows immediately from Lemma~\ref{4.1}.

Let $f : R \rightarrow \hat{R}$ be the completion.
By Lemma~4.1 in \cite{K16}, 
we have a map $f^* : \chow{*}(R)\subq \rightarrow \chow{*}(\hat{R})\subq$
that makes the following diagram commutative:
\bangou
\begin{equation}\label{eq3}
\begin{array}{ccc}
\g(R)\subq & \stackrel{\tau_{R/S}}{\ya{aaaa}} & \chow{*}(R)\subq \\
{\scriptstyle f^*}\downarrow\hphantom{\scriptstyle f^*} & & 
\hphantom{\scriptstyle f^*}\downarrow{\scriptstyle f^*} \\
\g(\hat{R})\subq & \stackrel{\tau_{\hat{R}/\hat{S}}}{\ya{aaaa}} & 
\chow{*}(\hat{R})\subq
\end{array}
\end{equation}
Note that the maps $f^* : \g(R)\subq \rightarrow \g(\hat{R})\subq$ and
$f^* : \chow{*}(R)\subq \rightarrow \chow{*}(\hat{R})\subq$ are
independent of the choice of $S$.
Furthermore, the map $\tau_{\hat{R}/\hat{S}}$ as above does not depend
on the choice of $\hat{S}$ 
since $\hat{R}$ is complete~(cf.\ Section~4 in \cite{K16}).

By Lemma~\ref{4.1}, it is easy to see that there is an induced map
$\ol{f^*} : \ol{\chow{*}(R)\subq} \rightarrow \ol{\chow{*}(\hat{R})\subq}$
that is independent of the choice of $S$, and
makes the following diagrams commutative:
\bangou
\begin{equation}\label{4.3}
\begin{array}{cccccccc}
\chow{*}(R)\subq & \stackrel{f^*}{\ya{aaaa}} & \chow{*}(\hat{R})\subq 
& & & 
\ol{\g(R)\subq} & \stackrel{\ol{\tau_{R/S}}}{\ya{aaaa}} & \ol{\chow{*}(R)\subq} \\
\downarrow & & \downarrow & & & 
{\scriptstyle \ol{f^*}}\downarrow\hphantom{\scriptstyle \ol{f^*}} & & 
\hphantom{\scriptstyle \ol{f^*}}\downarrow{\scriptstyle \ol{f^*}} \\
\ol{\chow{*}(R)\subq} & \stackrel{\ol{f^*}}{\ya{aaaa}} & \ol{\chow{*}(\hat{R})\subq} 
& & & 
\ol{\g(\hat{R})\subq} & \stackrel{\ol{\tau_{\hat{R}/\hat{S}}}}{\ya{aaaa}} & 
\ol{\chow{*}(\hat{R})\subq} 
\end{array}
\end{equation}
Since $\tau_{\hat{R}/\hat{S}}$ is independent of the choice of $S$,
so is $\ol{\tau_{\hat{R}/\hat{S}}}$. 
Since both $\ol{\tau_{R/S}}$ and $\ol{\tau_{\hat{R}/\hat{S}}}$ are 
isomorphisms (see Section~\ref{2}) and 
$\ol{f^*} : \ol{\g(R)\subq} \rightarrow \ol{\g(\hat{R})\subq}$
is injective by (1),
it follows that the map $\ol{\tau_{R/S}} : \ol{\g(R)\subq} \rightarrow
\ol{\chow{*}(R)\subq}$ is independent of the choice of $S$.
\qed

We hereafter denote the map $\ol{\tau_{R/S}}$ simply by $\ol{\tau_R}$.

\vspace{2mm}

Note that, if the maps $f^*$ in the diagram (\ref{eq3}) are injective,
then $\tau_{R/S}$ itself is independent of the choice of $S$.
The author does not know any example such that $f^*$ is not injective.
We refer the reader to \cite{KK} for some sufficient conditions 
of the injectivity of the map $f^* : \g(R)\subq \rightarrow \g(\hat{R})\subq$.

\vspace{2mm}

We proved that $\ol{f^*} : \ol{\g(R)\subq} \rightarrow \ol{\g(\hat{R})\subq}$
is injective in Theorem~\ref{4.2} (1).
The remainder of this section will be devoted to investigating some 
sufficient conditions of the surjectivity of the map
$\ol{f^*} : \ol{\g(R)\subq} \rightarrow \ol{\g(\hat{R})\subq}$.
Note that $\ol{f^*} : \ol{\g(R)\subq} \rightarrow \ol{\g(\hat{R})\subq}$
is surjective if and only if so is
$\ol{f^*} : \ol{\chow{*}(R)\subq} \rightarrow \ol{\chow{*}(\hat{R})\subq}$.

Here, we give an easy example where $\ol{f^*}$ is not surjective.
Set $R = {\Bbb C}[x,y]_{(x,y)}/(xy - x^3 - y^3)$.
Since $xy - x^3 - y^3$ is an irreducible polynomial, we have
\[
{\Bbb Q} = \chow{1}(R)\subq = \ol{\chow{1}(R)\subq}
\] 
by Proposition~\ref{3.5}~(3).
However, since the ring $\hat{R} = {\Bbb C}[[x,y]]/(xy - x^3 - y^3)$
has two minimal prime ideals, we have
\[
{\Bbb Q}\oplus{\Bbb Q} = \chow{1}(\hat{R})\subq = \ol{\chow{1}(\hat{R})\subq} 
\] 
by Proposition~\ref{3.5}~(3).

\bangou
\begin{Remark}\label{4.4}
\begin{rm}
Let $(R,\m)$ be a Noetherian local ring that satisfies the assumptions
in Theorem~\ref{3.1}.

Then, $\ol{f^*} : \ol{\g(R)\subq} \rightarrow \ol{\g(\hat{R})\subq}$ is an isomorphism
if and only if the equality 
$\dim \ol{\k{\m}(R)\subq} = \dim \ol{\k{\m\hat{R}}(\hat{R})\subq}$ holds.

On the other hand, 
the natural map $f^* : \k{\m}(R)\subq \rightarrow 
\k{\m\hat{R}}(\hat{R})\subq$ is an 
isomorphism by Theorem~7.1 in Thomason-Trobaugh~\cite{TT}.
By the surjectivity of $f^*$,
it is easy to see that $f^*(\N\k{\m}(R)\subq) \supseteq 
\N\k{\m\hat{R}}(\hat{R})\subq$.

Therefore, we know that 
$\ol{f^*} : \ol{\g(R)\subq} \rightarrow \ol{\g(\hat{R})\subq}$ 
is an isomorphism if and only if 
$f^*(\N\k{\m}(R)\subq) \subseteq \N\k{\m\hat{R}}(\hat{R})\subq$.
\end{rm}
\end{Remark}

\bangou
\begin{Proposition}\label{4.5}
Assume that $R$ is a local ring that satisfies the assumptions in
Theorem~\ref{3.1}.
Furthermore, assume that $R$ is henselian or that $R$ is the local ring 
(at the homogeneous maximal ideal) of an affine cone 
of a smooth projective variety over a field.
Then, $\ol{f^*} : \ol{\g(R)\subq} \rightarrow \ol{\g(\hat{R})\subq}$ 
is an isomorphism.
\end{Proposition}

\proof
First assume that $R$ is henselian.
Let ${\Bbb F}.$ be a complex contained in $C^{\m}(R)$.
Let $M$ be a finitely generated $\hat{R}$-module.
Then, using the approximation theorem due to Popescu-Ogoma~(\cite{ogoma},
\cite{P}),
there exists a finitely generated $R$-module $N$ 
such that $\chi_{{\Bbb F}.}([N]) = \chi_{{\Bbb F}.\otimes_R\hat{R}}([M])$.
Here, assume that ${\Bbb F}. \in \N\k{\m}(R)\subq$.
Then, we have $\chi_{{\Bbb F}.}([N]) = 0$.
Therefore, we have ${\Bbb F}.\otimes_R\hat{R} \in 
\N\k{\m\hat{R}}(\hat{R})\subq$.

Next, assume that $R$ is the local ring at the homogeneous maximal ideal
of an affine cone of a smooth projective variety over a field, i.e.,
let $A = \oplus_{n \geq 0}A_n = k[A_1]$ be a standard Noetherian graded ring 
over a field $k$ such that $X = \proj(A)$ is smooth over $k$, and
set $R = A_{A_+}$, $\m = A_+R$, where $A_+ = \oplus_{n > 0}A_n$.
Let $\pi : Z \rightarrow \spec(R)$ be the blowing-up with center $\m$.
Then, $\pi$ gives a resolution of singularities of $\spec(R)$.
Since $R \rightarrow \hat{R}$ is a regular homomorphism,
$\pi\times 1 : Z\times_{\spec(R)}\spec(\hat{R}) \rightarrow \spec(\hat{R})$ 
is also a resolution of singularities.
It is easy to see that $X$ is isomorphic to $\pi^{-1}(\m)$.
We have the following two fibre squares:
\[
\begin{array}{cccccccc}
X & \stackrel{j}{\rightarrow} & Z & & & X & \rightarrow & 
Z\times_{\spec(R)}\spec(\hat{R}) \\
{\scriptstyle \pi'}\downarrow\hphantom{\scriptstyle \pi'} & & 
{\scriptstyle \pi}\downarrow\hphantom{\scriptstyle \pi} & & & 
\downarrow & & {\scriptstyle \pi\times 1}\downarrow\hphantom{\scriptstyle 
\pi\times 1} \\
\spec(R/\m) & \rightarrow & \spec(R) & & & 
\spec(\hat{R}/\m\hat{R}) & \rightarrow
& \spec(\hat{R})
\end{array}
\]
By the argument in Roberts-Srinivas~\cite{RS} (see (\ref{6.3}) below),
for $\alpha \in \k{\m}(R)\subq$, $\alpha$ is contained in $\N\k{\m}(R)\subq$
if and only if $\alpha$ is contained in the kernel of the composite map 
\[
\k{\m}(R)\subq \stackrel{\pi^*}{\longrightarrow} \k{X}(Z)\subq
\stackrel{\chi}{\longrightarrow} \g(X)\subq
\stackrel{\tau_{X/Z}}{\longrightarrow} {\rm CH}^\cdot(X)\subq
\longrightarrow {\rm CH}^\cdot_{\rm num}(X)\subq ,
\]
where $\chi$ is a map defined in the proof of Claim~\ref{3.3},
and $\tau_{X/Z}$ is an isomorphism given by the singular Riemann-Roch theorem
(Chapter~18 in \cite{F}) with base regular scheme $Z$.

In the same way, one can prove that 
$f^*(\alpha)$ is contained in $\N\k{\m\hat{R}}(\hat{R})\subq$
if and only if $f^*(\alpha)$ is contained in the kernel 
of the composite map 
\[
\k{\m\hat{R}}(\hat{R})\subq \stackrel{(\pi\times 1)^*}{\ya{aaaaa}} 
\k{X}(Z\times \spec(\hat{R}))\subq
\stackrel{\chi}{\longrightarrow} \g(X)\subq
\stackrel{\tau_{X/(Z\times \spec(\hat{R}))}}{\ya{aaaaaaaaaa}} 
{\rm CH}^\cdot(X)\subq
\longrightarrow {\rm CH}^\cdot_{\rm num}(X)\subq .
\]
By Corollary~18.1.2 in \cite{F} and the definition of the maps, we have
\[
\tau_{X/Z} = {\rm td}(N)^{-1} \cdot \tau_{X/X}
= \tau_{X/(Z\times \spec(\hat{R}))} ,
\]
where $N$ is the normal bundle of the embedding
$X \rightarrow Z$.
Note that the diagram
\[
\begin{array}{ccccc}
\k{\m}(R)\subq & \stackrel{\pi^*}{\longrightarrow} & \k{X}(Z)\subq & & \\
{\scriptstyle f^*}\downarrow\hphantom{\scriptstyle f^*} & & 
{\scriptstyle (1 \times f)^*}\downarrow\hphantom{\scriptstyle 
(1 \times f)^*} & 
\searrow & \\
\k{\m\hat{R}}(\hat{R})\subq & \stackrel{(\pi\times 1)^*}{\longrightarrow} &
\k{X}(Z\times \spec(\hat{R}))\subq & \longrightarrow & \g(X)\subq
\end{array}
\]
is commutative.
Thus, if $\alpha$ is contained in $\N\k{\m}(R)\subq$, then
$f^*(\alpha)$ is in $\N\k{\m\hat{R}}(\hat{R})\subq$.
\qed

\section{Numerically Roberts rings}\label{5}

Since $(R,\m)$ is a homomorphic image of a regular local ring $S$,
we have an isomorphism
\[
\ol{\tau_R} : \ol{\g(R)\subq} \rightarrow \ol{\chow{*}(R)\subq}
= \bigoplus_{i = 0}^{\dim R}\ol{\chow{i}(R)\subq}
\]
that is independent of the choice of a base regular local ring $S$,
as was shown in Section~\ref{4}.
Using the map as above, we shall define a notion of {\em numerically Roberts}
rings and study these rings in this section.

\bangou
\begin{Definition}\label{5.1}
\begin{rm}
We say that $R$ is a {\em numerically Roberts} ring if
$\ol{\tau_R}(\ol{[R]}) \in \ol{\chow{\dim R}(R)\subq}$.
\end{rm}
\end{Definition}

Set $d = \dim R$ and 
\bangou
\begin{equation}\label{5.2}
\tau_{R/S}([R]) = \tau_d + \cdots + \tau_0 ,
\end{equation}
where $\tau_i \in \chow{i}(R)\subq$ for $i = 0, \cdots, d$.
By the top term property (Theorem~18.3~(5) in \cite{F}), we have
\[
\tau_d = \sum_{\p \in \assh(R)} \ell_{R_\p}(R_\p)[\spec(R/\p)] .
\]
As in Definition~2.1 in \cite{K16}, 
$R$ is said to be a {\em Roberts} ring if
$\tau_{d-1} = \cdots = \tau_0 = 0$ with some base regular local ring $S$.
Since the diagram
\[
\begin{array}{ccccc}
\g(R)\subq & \stackrel{\tau_{R/S}}{\ya{aaaa}} & \chow{*}(R)\subq & = &
\oplus_{i = 0}^d \chow{i}(R)\subq \\
\downarrow & & \downarrow & & \\
\ol{\g(R)\subq} & \stackrel{\ol{\tau_R}}{\ya{aaaa}} & \ol{\chow{*}(R)\subq} 
& = & \oplus_{i = 0}^d \ol{\chow{i}(R)\subq}
\end{array}
\]
is commutative, we have 
\bangou
\begin{equation}\label{5.3}
\ol{\tau_R}(\ol{[R]}) = \ol{\tau_d} + \cdots + \ol{\tau_0} .
\end{equation}
Therefore, $R$ is a numerically Roberts ring if and only if 
$\tau_i \in \N\chow{i}(R)\subq$ for $i = 0, 1, \ldots, d-1$.
In particular, if $R$ is a Roberts ring, then 
it is a numerically Roberts ring.
(The converse is not true. See Example~\ref{5.6}.)

Assume that the natural map $\chow{*}(R)\subq \rightarrow \ol{\chow{*}(R)\subq}$
is an isomorphism.
Then, $R$ is a numerically Roberts ring if and only if 
$R$ is a Roberts ring.

For a complex ${\Bbb F}. \in C^\m(R)$,
the rational number 
\[
{\rm ch}({\Bbb F}.)(\tau_d) = {\rm ch}_d({\Bbb F}.)(\tau_d)
\]
is called the {\em Dutta multiplicity} of the complex ${\Bbb F}.$
and is denoted by $\chi_\infty({\Bbb F}.)$.

The following proposition characterizes numerically Roberts rings.
We refer the reader to \cite{Mo} for Hilbert-Kunz multiplicity.

\bangou
\begin{Theorem}\label{5.4}
\DuttaHilbertKunz
\end{Theorem}

\proof
With notation as in (\ref{5.2}), 
$R$ is a numerically Roberts ring if and only if
$(\tau_{R/S})^{-1}(\tau_{d-1} + \cdots + \tau_0) \in \N\g(R)\subq$
(cf.\ Proposition~\ref{2.3}).
Note that
\begin{eqnarray*}
\chi_\infty({\Bbb F}.) & = & 
\chi_{{\Bbb F}.} \left( (\tau_{R/S})^{-1}(\tau_d) \right) \\
& = & \chi_{{\Bbb F}.} \left( 
[R] - (\tau_{R/S})^{-1}(\tau_{d-1} + \cdots + \tau_0) \right) \\
& = & \chi({\Bbb F}.) - \chi_{{\Bbb F}.} \left( 
(\tau_{R/S})^{-1}(\tau_{d-1} + \cdots + \tau_0) \right) .
\end{eqnarray*}
Therefore, $R$ is a numerically Roberts ring if and only if 
$\chi_\infty({\Bbb F}.)$ coincides with $\chi({\Bbb F}.)$
for any ${\Bbb F}. \in C^\m(R)$.

Assume that $R$ is a Cohen-Macaulay ring of characteristic $p$.
Let $F(R)$ be the category of $R$-modules of finite length and
of finite projective dimension.
Since $R$ is a Cohen-Macaulay local ring, 
$\k{\m}(R)\subq$ is generated by free resolutions
of modules in $F(R)$ (cf.\ Proposition~2 in \cite{RS}).
For $M \in F(R)$, ${{\Bbb F}_M}. \in C^\m(R)$ denotes 
the minimal free resolution of $M$.
By (1), $R$ is a numerically Roberts ring if and only if 
$\chi_\infty({{\Bbb F}_M}.)$ coincides with $\chi({{\Bbb F}_M}.) = \ell_R(M)$
for any $M \in F(R)$.
Suppose that $M \in F(R)$.
Using the method of Lemma~9.10 in \cite{Sr}, 
$\m$-primary ideals $J$, $I_1$, \ldots, $I_t$ of finite projective dimension
can be found 
such that $I_1$, \ldots, $I_t$ are parameter ideals and 
\[
[M] = [R/J] - \sum_{i = 1}^t[R/I_i]
\]
in ${\rm K}_0(F(R))$.
Then, we have
\[
\chi([{{\Bbb F}_M .}]) = \ell_R(M) = \ell_R(R/J) - \sum_{i = 1}^t \ell_R(R/I_i)
\]
and
\[
[{{\Bbb F}_M .}] = 
[{{\Bbb F}_{R/J}}.] - \sum_{i = 1}^t[{{\Bbb F}_{R/I_i}}.]
\]
in $\k{\m}(R)\subq$.
Therefore, we have
\begin{eqnarray*}
\chi_\infty({{\Bbb F}_M}.) & = & 
\chi_\infty({{\Bbb F}_{R/J}}.)-\sum_{i = 1}^t\chi_\infty({{\Bbb F}_{R/I_i}}.)\\
& = & e_{HK}(J)-\sum_{i = 1}^t\ell_R(R/I_i) ,
\end{eqnarray*}
by Remark~2.7 in \cite{K9} and Theorem~1.2~(1) in \cite{K15}.
Thus, $R$ is a numerically Roberts ring if and only if 
$e_{HK}(J)$ coincides with $\ell_R(R/J)$ for any $\m$-primary ideal $J$
of finite projective dimension.
\qed

\bangou
\begin{Remark}\label{newrem6.5}
\begin{rm}
Let $R$ be a numerically Roberts ring.
Then, $R$ satisfies the vanishing property of intersection multiplicities,
that is,
\[
\sum_i(-1)^i\ell_R(\tor{i}{R}{M}{N}) = 0
\]
for finitely generated $R$-modules $M$ and $N$ 
such that 
$\hd_RM < \infty$, $\hd_RN < \infty$, 
$\ell_R(M\otimes_RN) < \infty$ and $\dim M + \dim N < \dim R$.

The proof is the same as that due to Roberts~\cite{R1}.
\end{rm}
\end{Remark}

\bangou
\begin{Example}\label{5.5}
\begin{rm}
Let $R$ be a homomorphic image of a regular local ring and let $d = \dim R$.

Suppose $d = 0$.
Since $\ol{\chow{*}(R)\subq} = \ol{\chow{0}(R)\subq}$,
$R$ is a numerically Roberts ring.

Suppose $d = 1$.
Since $\ol{\chow{0}(R)\subq} = 0$ by Proposition~\ref{3.5}~(1),
$R$ is a numerically Roberts ring.

Assume that $R$ is equi-dimensional with $d = 2$.
Since $\ol{\chow{0}(R)\subq} = \ol{\chow{1}(R)\subq} = 0$ by 
Proposition~\ref{3.5},
$R$ is a numerically Roberts ring.

Assume that $R$ is a Gorenstein ring with $d = 3$.
Set
\[
\ol{\tau_R}(\ol{[R]}) = \ol{\tau_3} + \ol{\tau_2} + \ol{\tau_1} + \ol{\tau_0} 
\]
as (\ref{5.3}).
Since $\ol{\chow{0}(R)\subq} = \ol{\chow{1}(R)\subq} = 0$ by Proposition~\ref{3.5},
we have $\ol{\tau_0} = \ol{\tau_1} = 0$.
Furthermore, since $R$ is a Gorenstein ring, we have $\tau_2 = 0$
by Proposition~2.8 in \cite{K11}.
Therefore, $R$ is a numerically Roberts ring.

Using an example due to Dutta-Hochster-MacLaughlin~\cite{DHM},
we can construct a three-dimensional Cohen-Macaulay normal 
ring $R$ and ${\Bbb F}. 
\in C^\m(R)$ such that $\chi_\infty({\Bbb F}.) \neq \chi({\Bbb F}.)$.
By Proposition~\ref{5.4} (1), $R$ is not a numerically Roberts ring.

A five-dimensional Gorenstein ring constructed by Miller-Singh~\cite{MS}
is not a numerically Roberts ring.
\end{rm}
\end{Example}

Next, we give an example of a numerically Roberts ring 
that is not a Roberts ring.

\bangou
\begin{Example}\label{5.6}
\begin{rm}
We give an example of a two-dimensional Noetherian local domain
that is not a Roberts ring.
(Recall that two-dimensional Noetherian local domains
are numerically Roberts rings by Example~\ref{5.5}.)

Let $X$ be a smooth projective curve over ${\Bbb C}$
with $g(X) \geq 2$.
Since ${\rm Pic}^0(X)$ is an abelian variety of dimension $g(X)$,
${\rm Pic}(X)\subq$ is an infinite dimensional ${\Bbb Q}$-vector space.
Furthermore, we have $\deg K_X = 2 g(X) -2 > 0$.

Take a divisor $H \in {\rm Div}(X)$ such that
$H$ and $K_X$ are linearly independent in ${\rm Pic}(X)\subq$.
We assume $\deg H >>0$.
Then, $H$ is an very ample divisor,
and we may assume that the graded ring
\[
A = \bigoplus_{n \geq 0} H^0(X, {\cal O}_X(nH))
\]
is generated by elements of degree $1$ over ${\Bbb C} = 
H^0(X, {\cal O}_X)$.

Set $R = A_{A_+}$.
Then, $R$ is a two-dimensional Noetherian local domain
and we have an exact sequence
\[
{\rm CH}^\cdot(X)\subq \stackrel{H}{\ya{aaaa}}
{\rm CH}^\cdot(X)\subq \stackrel{\xi}{\ya{aaaa}}
\chow{*}(R)\subq \stackrel{}{\ya{aaaa}} 0 ,
\]
where $\xi$ is defined by $\xi([\proj(A/P)]) = [\spec(R/PR)]$
for each homogeneous prime ideal $P$ of $A$ not equal to $A_+$.
The map ${\rm CH}^\cdot(X)\subq \stackrel{H}{\ya{aaaa}}
{\rm CH}^\cdot(X)\subq$ denotes the multiplication by $H$.
Then, by Theorem~1.3 in \cite{K11}, we have 
\[
\xi({\rm td}(\Omega_X^\vee)) = \tau_R([R]) = \tau_2 + \tau_1 + \tau_0 .
\]
Furthermore, by definition of Todd classes of vector bundles,
we have
\[
{\rm td}(\Omega_X^\vee) = 1 + \frac{1}{2}c_1(\Omega_X^\vee) =
1 - \frac{1}{2}c_1(\omega_X) = 1 - \frac{1}{2}K_X .
\]
Therefore, we have $\xi(- \frac{1}{2}K_X) = \tau_1$.
Since $H$ and $K_X$ are linearly independent in ${\rm CH}^1(X)\subq$,
we obtain $\tau_1 \neq 0$.
Hence $R$ is not a Roberts ring.
\end{rm}
\end{Example}

\bangou
\begin{Remark}\label{5.7}
\begin{rm}
The author knows no example of a four-dimensional Gorenstein ring 
that is not a numerically Roberts ring.

Let $R$ be a four-dimensional Gorenstein ring and set 
\[
\ol{\tau_R}(\ol{[R]}) = \ol{\tau_4} + \ol{\tau_3} +  \ol{\tau_2} + 
\ol{\tau_1} + \ol{\tau_0} .
\]
Since $\ol{\chow{0}(R)\subq} = \ol{\chow{1}(R)\subq} = 0$ by Proposition~\ref{3.5},
we have $\ol{\tau_0} = \ol{\tau_1} = 0$.
Since $R$ is a Gorenstein ring, we have $\tau_3 = 0$.
Therefore, $R$ is a numerically Roberts ring if and only if $\ol{\tau_2} = 0$.

Furthermore, assume that $R$ is the local ring 
(at the homogeneous maximal ideal) of 
an affine cone  of a smooth projective 
variety $X$ of dimension $3$ over ${\Bbb C}$.
In this case, it will be proved in Remark~\ref{6.9} that
$\ol{\chow{2}(R)\subq} = 0$.
Therefore, 
$R$ is a numerically Roberts ring in this case.
\end{rm}
\end{Remark}


\bangou
\begin{Remark}\label{5.9}
\begin{rm}
Assume that $R$ is a numerically Roberts ring.
Then, $R/xR$ is also a numerically Roberts ring
for any non-zero-divisor $x$ of $R$, as will be proved in this remark.
The author does not know whether $R_\p$ is a numerically Roberts ring
for a prime ideal $\p$ of $R$.

Let $(R,\m)$ be a homomorphic image of a regular local ring $S$.
Assume that $(A,\n)$ is a homomorphic image of $R$ such that $\hd_RA < \infty$.
Let ${\Bbb H}.$ be a finite $R$-free resolution of $A$.
Then, we have the map
\[
\chi_{{\Bbb H}.} : \g(R)\subq \rightarrow \g(A)\subq
\]
defined by $\chi_{{\Bbb H}.}([M]) = \sum_i(-1)^i[H_i({\Bbb H}.
\otimes_RM)]$.
Suppose that ${\Bbb F}.$ is a complex in $C^\n(A)$.
Then, since $\hd_RA < \infty$, there is a complex ${\Bbb G}. \in C^\m(R)$ 
with a quasi-isomorphism ${\Bbb G}. \rightarrow {\Bbb F}.$ as in 
Lemma~1.10 in \cite{GS}.
Using a spectral sequence argument,
the diagram
\bangou
\begin{equation}\label{5.10}
\begin{array}{ccc}
\g(R)\subq & \stackrel{\chi_{{\Bbb H}.}}{\ya{aaaaa}} & \g(A)\subq \\
& {\scriptstyle \chi_{{\Bbb G}.}}\searrow\hphantom{\scriptstyle 
\chi_{{\Bbb G}.}} & 
\hphantom{\scriptstyle \chi_{{\Bbb F}.}}\downarrow{\scriptstyle 
\chi_{{\Bbb F}.}} \\
& & {\Bbb Q}
\end{array}
\end{equation}
is commutative.
Therefore, for $\beta \in \N\g(R)\subq$, $\chi_{{\Bbb H}.}(\beta)$ 
is contained in $\N\g(A)\subq$.
Thus, we have an induced map
$\ol{\chi_{{\Bbb H}.}} : \ol{\g(R)\subq} \rightarrow \ol{\g(A)\subq}$ 
that makes the following diagram commutative:
\[
\begin{array}{ccc}
\g(R)\subq & \stackrel{\chi_{{\Bbb H}.}}{\ya{aaaaa}} & \g(A)\subq \\
\downarrow & & \downarrow \\
\ol{\g(R)\subq} & \stackrel{\ol{\chi_{{\Bbb H}.}}}{\ya{aaaaa}} & \ol{\g(A)\subq} \\
\end{array}
\]
Since the diagram~(\ref{5.10}) is commutative,
so is
\[
\begin{array}{ccc}
\chow{*}(R)\subq & \stackrel{{\rm ch}({\Bbb H}.)}{\ya{aaaaa}} & \chow{*}(A)\subq \\
& 
\searrow
& 
\hphantom{\scriptstyle {\rm ch}({\Bbb F}.)}\downarrow{\scriptstyle 
{\rm ch}({\Bbb F}.)} \\
& {\scriptstyle {\rm ch}({\Bbb G}.)} & {\Bbb Q}
\end{array}
\]
by the local Riemann-Roch formula (Example~18.3.12 in \cite{F}).
Therefore, we have ${\rm ch}({\Bbb H}.)(\N\chow{*}(R)\subq) \subseteq 
\N\chow{*}(A)\subq$ and the following commutative diagram:
\bangou
\begin{equation}\label{5.11}
\begin{array}{ccc}
\chow{*}(R)\subq & \stackrel{{\rm ch}({\Bbb H}.)}{\ya{aaaaa}} & \chow{*}(A)\subq \\
\downarrow & & \downarrow \\
\ol{\chow{*}(R)\subq} & \stackrel{\ol{{\rm ch}({\Bbb H}.)}}{\ya{aaaaa}} & 
\ol{\chow{*}(A)\subq} \\
\end{array}
\end{equation}
Then, the following diagram is commutative:
\[
\begin{array}{ccc}
\ol{\g(R)\subq} & \stackrel{\ol{\tau_{R}}}{\ya{aaaaa}} & \ol{\chow{*}(R)\subq} \\
{\scriptstyle \ol{\chi_{{\Bbb H}.}}}\downarrow\hphantom{\scriptstyle 
\chi_{{\Bbb H}.}} & & \hphantom{\scriptstyle 
{\rm ch}({\Bbb H}.)}\downarrow{\scriptstyle \ol{{\rm ch}({\Bbb H}.)}} \\
\ol{\g(A)\subq} & \stackrel{\ol{\tau_{A}}}{\ya{aaaaa}} & \ol{\chow{*}(A)\subq} \\
\end{array}
\]
Therefore, we have 
\bangou
\begin{equation}\label{5.12}
\ol{{\rm ch}({\Bbb H}.)}\left( \ol{\tau_{R}}(\ol{[R]}) \right)
= \ol{\tau_{A}}(\ol{[A]}) .
\end{equation}

Now, assume that $A$ coincides with $R/xR$ for some non-zero-divisor
$x$ on $R$.
Then, for $\gamma \in \chow{i}(R)\subq$, ${\rm ch}({\Bbb H}.)(\gamma)$
is contained in $\chow{i-1}(A)\subq$ since 
${\rm ch}({\Bbb H}.) = {\rm ch}_1({\Bbb H}.)$
by Corollary~18.1.2 in Fulton~\cite{F}.
By the commutativity of the diagram~(\ref{5.11}), we have
\[
\ol{{\rm ch}({\Bbb H}.)}\left( \ol{\chow{i}(R)\subq} \right)
\subseteq \ol{\chow{i-1}(A)\subq}
\]
for each $i$.
Here, assume that $R$ is a $d$-dimensional numerically Roberts ring.
Since $\ol{\tau_{R}}(\ol{[R]}) \in \ol{\chow{d}(R)\subq}$, we have
\[
\ol{\tau_{R}}(\ol{[A]}) = \ol{{\rm ch}({\Bbb H}.)}\left( 
\ol{\tau_{R}}(\ol{[R]}) \right) \in \ol{\chow{d-1}(A)\subq}
\]
by (\ref{5.12}). 
Consequently $A$ is a $(d-1)$-dimensional numerically Roberts ring.
\end{rm}
\end{Remark}

\bangou
\begin{Remark}\label{5.13}
\begin{rm}
Let $f : (A,\p) \rightarrow (B,\q)$ be a flat local homomorphism 
of Noetherian local rings.
Assume that $B/\p B$ is a complete intersection, and 
$B/\q$ is a finite (algebraic) separable extension of $A/\p$.
Then, $A$ is a numerically Roberts ring
if and only if $B$ is a numerically Roberts ring.

The proof is omitted here.
\end{rm}
\end{Remark}

\section{Affine cones of smooth projective varieties}\label{6}
In this section, we treat affine cones of smooth projective 
varieties.
Using a method and a result of Roberts and Srinivas~\cite{RS},
we shall calculate Chow group modulo numerical equivalence.

Let $k$ be an algebraically closed field and let 
$A = \oplus_{n \geq 0}A_n$ be a Noetherian graded ring with
$A_0 = k$ and $A = k[A_1]$.
Assume that $X = \proj(A)$ is smooth over $k$.
Set $A_+ = \oplus_{n > 0}A_n$, $R = A_{A_+}$ and $\m = A_+R$.
Let $\pi : Z \rightarrow \spec(R)$ be the blowing-up with center $\m$.
Then, $\pi^{-1}(\m)$ naturally coincides with $X$.
Thus, we regard $X$ as a closed subscheme of $Z$.
Let $h$ be the very ample divisor under the embedding $X = \proj(A)$.

Note that $Z \setminus X = \spec(R) \setminus \spec(R/\m)$.
Then, by the theory of localization sequences 
due to Thomason-Trobaugh~\cite{TT},
we have the following commutative diagram
\[
\begin{array}{ccccccccc}
\cdots & \longrightarrow & {\rm K}_1(Z \setminus X)\subq & \longrightarrow &
\k{X}(Z)\subq & \stackrel{r}{\longrightarrow} & 
\k{}(Z)\subq & \longrightarrow & \cdots \\
& & \| & & {\scriptstyle \pi^*} \uparrow \hphantom{\scriptstyle \pi^*}
& & \uparrow & & \\
\cdots & \longrightarrow & {\rm K}_1(\spec(R) \setminus \spec(R/\m))\subq 
& \longrightarrow &
\k{\m}(R)\subq & \stackrel{s}{\longrightarrow} & 
\k{}(R)\subq & \longrightarrow & \cdots 
\end{array}
\]
where the horizontal sequences are exact, and
$\k{}(R)\subq$ (resp.\ $\k{}(Z)\subq$) denotes the Grothendieck group
of finitely generated projective $R$-modules
(resp.\ locally free ${\cal O}_Z$-modules of finite rank).
We refer the reader to 1.5 in \cite{GS} for the definition of $\pi^*$.
Since the map $s$ coincides with $0$ in this case, we have an exact sequence
\[
\k{\m}(R)\subq \stackrel{\pi^*}{\longrightarrow}
\k{X}(Z)\subq \stackrel{r}{\longrightarrow}
\k{}(Z)\subq .
\]

Since $Z$ is a regular scheme,
$\chi : \k{X}(Z)\subq \rightarrow \g(X)\subq$ is an isomorphism by 
Lemma~1.9 in \cite{GS},
where $\chi$ is defined by $\chi([{\Bbb F}.])
= \sum_i(-1)^i[H_i({\Bbb F}.)]$.
In particular, the natural map $\k{}(Z)\subq \rightarrow \g(Z)\subq$
is an isomorphism.
Therefore, we have an exact sequence
\[
\k{\m}(R)\subq \stackrel{\chi \pi^*}{\ya{aaa\chi \pi^*}}
\g(X)\subq \stackrel{i_*}{\longrightarrow}
\g(Z)\subq ,
\]
where $i_*$ is the map induced by the closed immersion 
$X \stackrel{i}{\longrightarrow} Z$.
Note that the diagram
\[
\begin{array}{ccc}
\g(X)\subq & \stackrel{i_*}{\longrightarrow} & \g(Z)\subq \\
{\scriptstyle \tau_{X/Z}} \downarrow 
\hphantom{{\scriptstyle \tau_{X/Z}}}
& & 
{\scriptstyle \tau_{Z/Z}} \downarrow 
\hphantom{{\scriptstyle \tau_{Z/Z}}} \\
\chow{*}(X)\subq & \stackrel{i_*}{\longrightarrow} & \chow{*}(Z)\subq
\end{array}
\]
is commutative, where $\tau_{X/Z}$ is the Riemann-Roch map of $X$
with base regular scheme $Z$ (see 18.2 and 20.1 in \cite{F}).
Since the vertical maps in the diagram as above are isomorphisms, the sequence
\[
\k{\m}(R)\subq \stackrel{\tau_{X/Z} \chi \pi^*}{\ya{aaa\tau_{X/Z} \chi \pi^*}}
\chow{*}(X)\subq \stackrel{i_*}{\longrightarrow}
\chow{*}(Z)\subq 
\]
is exact.

Let $Y$ be the blowing-up of $\spec(A)$ with center $A_+$.
Consider the following fibre squares:
\[
\begin{array}{ccccc}
X & \stackrel{i}{\longrightarrow} & Z & \stackrel{j}{\longrightarrow} & Y \\
{\scriptstyle \pi'} \downarrow 
\hphantom{{\scriptstyle \pi'}}
& & 
{\scriptstyle \pi} \downarrow 
\hphantom{{\scriptstyle \pi}}
& & 
{\scriptstyle \pi"} \downarrow 
\hphantom{{\scriptstyle \pi"}} \\
\spec(R/\m) & \longrightarrow & \spec(R) & \longrightarrow & \spec(A)
\end{array}
\]
Note that
$Y$ is a vector bundle on $X$ with sheaf of sections ${\cal O}_X(-1)$.
Let $p: Y \rightarrow X$ be the projection.
By Theorem~3.3~(a) in \cite{F}, the pull-back map 
$p^* : \chow{*}(X) \rightarrow \chow{*}(Y)$ is an isomorphism.
Consider the following commutative diagram:
\bangou
\begin{equation}\label{6.1}
\begin{array}{ccc}
\chow{*}(X)\subq & \stackrel{i_*}{\longrightarrow} & \chow{*}(Z)\subq \\
{\scriptstyle (ji)_*} \downarrow 
\hphantom{{\scriptstyle (ji)_*}}
& {\scriptstyle j^*} \nearrow 
\hphantom{{\scriptstyle j^*}} & 
{\scriptstyle i^*} \downarrow 
\hphantom{{\scriptstyle i^*}} \\
\chow{*}(Y)\subq & \stackrel{(ji)^*}{\longrightarrow} & \chow{*}(X)\subq
\end{array}
\end{equation}
Here, $i_*$ (resp.\ $(ji)_*$) is the map induced by the closed immersion
$i : X \rightarrow Z$ (resp.\ $ji : X \rightarrow Y$),
$j^*$ is the map induced by the flat map $j : Z \rightarrow Y$, and
$i^*$ (resp.\ $(ji)^*$) is the map taking the intersection with 
the Cartier divisor $i : X \rightarrow Z$ (resp.\ $ji : X \rightarrow Y$).
We refer to Chapters~1 and 2 in \cite{F} for these induced maps as above.
Commutativity of the diagram~(\ref{6.1}) follows from Proposition~1.7 and 
Proposition~2.3~(d) in \cite{F}.
By Corollary~6.5 in \cite{F}, $(ji)^*$ coincides with ${p^*}^{-1}$.
In particular, $(ji)^*$ is an isomorphism.
By the commutativity of the diagram~(\ref{6.1}), we have
\[
\ker\left( \chow{*}(X)\subq \stackrel{i_*}{\rightarrow} 
\chow{*}(Z)\subq \right)
=
\ker\left( \chow{*}(X)\subq \stackrel{(ji)^*(ji)_*}{\ya{aaaaaaaaa}} 
\chow{*}(X)\subq \right) .
\]
Using Example~3.3.2 in \cite{F}, the map $(ji)^*(ji)_*$ coincides with
the multiplication by $-h \in \chow{\dim X - 1}(X)\subq$.
Let $\chring{X} = \oplus_{i = 0}^{\dim X}{\rm CH}^i(X)\subq$ be the Chow ring
of $X$ with rational coefficients.
Note that ${\rm CH}^i(X) = \chow{\dim X - i}(X)$ for each $i$.
Setting $\varphi = \tau_{X/Z} \chi \pi^*$,
we obtain an exact sequence
\bangou
\begin{equation}\label{6.2}
\k{\m}(R)\subq \stackrel{\varphi}{\longrightarrow}
\chring{X} \stackrel{h}{\longrightarrow}
\chring{X} ,
\end{equation}
where $\chring{X} \stackrel{h}{\longrightarrow} \chring{X}$ denotes the
multiplication by $h$.

Set
\begin{eqnarray*}
K & = & 
\kernel({\rm CH}^\cdot(X)\subq \stackrel{h}{\longrightarrow} {\rm CH}^\cdot(X)\subq)
\\
L & = &
\kernel(\chnum{X} \stackrel{h}{\longrightarrow} \chnum{X}) 
\\
M & = & 
\chnum{X} \left/ h \cdot \chnum{X} \right.
\end{eqnarray*}
where $\chnum{X}$ denotes the Chow ring of $X$ (with rational coefficients)
modulo numerical equivalence.
Then, we have a commutative diagram
\[
\begin{array}{ccccccccccc}
& & \k{\m}(R)\subq & & & & & & & & \\
& & {\scriptstyle \varphi} \downarrow \hphantom{\scriptstyle \varphi}
& & & & & & & & \\
0 & \longrightarrow & K & \longrightarrow & {\rm CH}^\cdot(X)\subq 
& \stackrel{h}{\longrightarrow} & {\rm CH}^\cdot(X)\subq &  
\stackrel{\xi}{\longrightarrow} & \chow{*}(R)\subq & \longrightarrow & 0 \\
& & \downarrow & & \downarrow & & \downarrow & & \downarrow & & \\
0 & \longrightarrow & L & \longrightarrow & \chnum{X}
& \stackrel{h}{\longrightarrow} & \chnum{X} &  
\longrightarrow & M & \longrightarrow & 0
\end{array}
\]
where the horizontal sequences are exact (e.g., \cite{K11}), 
and $\xi$ is defined by
$\xi([\proj(A/P)]) = [\spec(R/PR)]$ for each homogeneous prime 
ideal $P$ of $A$ such that $P \neq A_+$.
Note that $\varphi : \k{\m}(R)\subq \rightarrow K$ is surjective 
by the exactness of the sequence~(\ref{6.2}).
Let $W$ be the image of the induced map $K \rightarrow L$.
We denote by 
\bangou
\begin{equation}\label{phi}
\phi : \k{\m}(R)\subq \rightarrow W
\end{equation}
 the induced surjection.
By the argument of Roberts-Srinivas~\cite{RS},
we know that
\bangou
\begin{eqnarray}
\label{6.3} \ & & 
\mbox{$W = \ol{\k{\m}(R)\subq}$, that is, $0 \rightarrow \N\k{\m}(R)\subq
\rightarrow \k{\m}(R)\subq \stackrel{\phi}{\rightarrow} W \rightarrow 0$ is exact,
and} \\
\bangou \ & & \nonumber \\
\label{6.4} & & 
\mbox{there exists a map 
$\ol{\xi} : \chnum{X} \rightarrow \ol{\chow{*}(R)\subq}$
such that the following } \\ \ & & 
\mbox{diagram is commutative, where the vertical maps are 
the projections.}
\nonumber
\end{eqnarray}
\[
\begin{array}{ccc}
{\rm CH}^\cdot(X)\subq &  
\stackrel{\xi}{\longrightarrow} & \chow{*}(R)\subq \\
\downarrow & & \downarrow \\
\chnum{X} & \stackrel{\ol{\xi}}{\longrightarrow} & \ol{\chow{*}(R)\subq}
\end{array} 
\]

Here, we give an outline of proofs of (\ref{6.3}) and (\ref{6.4}).
For a homogeneous prime ideal $P$ of $A$ not equal to $A_+$, 
$Z_P$ denotes the proper transform of $\spec(R/PR)$,
i.e., $Z_P$ is the closed integral subscheme of $Z$ such that 
$\pi(Z_P) = \spec(R/PR)$.
Note that 
the induced morphism $Z_P \rightarrow \spec(R/PR)$ is a birational surjection.
Consider the following fibre square:
\[
\begin{array}{ccc}
X & \stackrel{i}{\longrightarrow} & Z \\
{\scriptstyle \pi'} \downarrow 
\hphantom{{\scriptstyle \pi'}}
& & 
{\scriptstyle \pi} \downarrow 
\hphantom{{\scriptstyle \pi}} \\
\spec(R/\m) & \longrightarrow & \spec(R)
\end{array}
\]

For $\alpha \in \k{\m}(R)\subq$, we have
\[
{\rm ch}(\alpha)\left( [\spec(R/PR)] \right)
={\rm ch}(\alpha)\left( \pi_*([Z_P]) \right)
=\pi'_*\left( {\rm ch}(\pi^*\alpha)\left( [Z_P] \right) \right)
\]
by Definition~17.1 in \cite{F}.
Since $Z$ and $X$ are regular schemes, we have isomorphisms
\[
\k{X}(Z)\subq \stackrel{\chi}{\longrightarrow} \g(X)\subq
\stackrel{\chi}{\longleftarrow} \k{X}(X)\subq 
\]
by Lemma~1.9 in \cite{GS}.
Take $\epsilon \in \k{X}(X)\subq$ such that 
$\chi(\pi^*\alpha) = \chi(\epsilon)$.
Since the diagram
\[
\begin{array}{ccc}
\g(X)\subq & \stackrel{\tau_{X/X}}{\longrightarrow} & \chow{*}(X)\subq \\
{\scriptstyle \chi_\epsilon} \downarrow 
\hphantom{{\scriptstyle \chi_\epsilon}}
& & 
{\scriptstyle {\rm ch}(\epsilon)} \downarrow 
\hphantom{{\scriptstyle {\rm ch}(\epsilon)}} \\
\g(X)\subq & \stackrel{\tau_{X/X}}{\longrightarrow} & \chow{*}(X)\subq
\end{array}
\]
is commutative by the local Riemann-Roch formula~(Example~18.3.12 in \cite{F}),
we have 
\[
\tau_{X/X} \chi (\pi^*\alpha) = \tau_{X/X} \chi (\epsilon)
= \tau_{X/X} \chi_\epsilon([{\cal O}_X]) = 
{\rm ch}(\epsilon)\left( \tau_{X/X}([{\cal O}_X]) \right)
= {\rm ch}(\epsilon)\left( [X] \right) .
\]
Note that, by the definition of the Riemann-Roch map (18.2 and 20.1 in 
\cite{F}), we have $\tau_{X/X}([{\cal O}_X]) = X$.

On the other hand, by Corollary~18.1.2 in \cite{F}, we have 
\[
\pi'_*\left( {\rm ch}(\pi^*\alpha)\left( [Z_P] \right)\right)
= \pi'_*\left( {\rm ch}(\epsilon) \cdot {\rm td}(N)^{-1} \cdot 
i^*[Z_P] \right) ,
\]
where $N$ denotes the normal bundle of $i : X \rightarrow Z$,
and $i^* : \chow{*}(Z)\subq \rightarrow \chow{*}(X)\subq$ denotes 
the map taking the intersection with $X$.
Since $i^*[Z_P] = [\proj(A/P)]$ and 
${\rm ch}(\epsilon)\left( [X] \right) = \tau_{X/X}\chi (\pi^*\alpha)$
as above,
we have
\[
\pi'_*\left( {\rm ch}(\epsilon) \cdot {\rm td}(N)^{-1} \cdot 
i^*[Z_P] \right)
=
\pi'_*\left( \tau_{X/X}\chi (\pi^*\alpha) \cdot {\rm td}(N)^{-1} \cdot 
[\proj(A/P)] \right) .
\]
By Corollary~18.1.2 in \cite{F}
and the definition of the Riemann-Roch map~(18.2 and 20.1 in \cite{F}),
we have
\[
\tau_{X/Z} = {\rm td}(N)^{-1}
\cdot \tau_{X/X} .
\]
Hence, we have
\begin{eqnarray*}
\pi'_*\left( \tau_{X/X}\chi (\pi^*\alpha) \cdot {\rm td}(N)^{-1} \cdot 
[\proj(A/P)] \right)
& = &
\pi'_*\left( \tau_{X/Z}\chi (\pi^*\alpha) \cdot [\proj(A/P)] \right) \\
& = & \pi'_*\left( \varphi(\alpha) \cdot [\proj(A/P)] \right) 
\end{eqnarray*}
by the definition of $\varphi$.
Consequently, we have
\bangou
\begin{equation}\label{6.5}
{\rm ch}(\alpha)\left( [\spec(R/PR)] \right)
=
\pi'_*\left( \varphi(\alpha) \cdot [\proj(A/P)] \right) .
\end{equation}

Recall that $\chow{*}(R)\subq$ is generated by
\[
\{ [\spec(R/PR)] \mid \mbox{$P$ is a homogeneous prime ideal of $A$ 
not equal to $A_+$} \} 
\]
because $\xi$ is surjective (cf.\ Theorem~1.3 in \cite{K11}).
Therefore, by (\ref{6.5}), 
\begin{quote}
$\alpha$ is in $\N\k{\m}(R)\subq$ 
if and only if $\varphi(\alpha)$ coincides with $0$
in $\chnum{X}$.
\end{quote}
In other words, $\N\k{\m}(R)\subq = \ker(\phi)$, where $\phi$
is the map in (\ref{phi}).
The statement (\ref{6.3}) follows from this equality.

If $[\proj(A/P)]$ is equal to $0$ in $\chnum{X}$,
then $[\spec(R/PR)]$ is contained in $\N\chow{*}(R)\subq$ by (\ref{6.5}).
Thus, (\ref{6.4}) holds.
\qed

By (\ref{6.4}) above, 
we have a surjection $M \rightarrow \ol{\chow{*}(R)\subq}$.
Since $\dim \chnum{X} < \infty$, we have $\dim L = \dim M$.
By (\ref{6.3}) as above and Theorem~\ref{3.1}, 
we have $\dim W = \dim \ol{\chow{*}(R)\subq}$.
\[
\begin{array}{ccc}
\dim W & \leq & \dim L \\
\| & & \| \\
\dim \ol{\chow{*}(R)\subq} & \leq & \dim M
\end{array}
\]
Therefore, the following three conditions are equivalent:
\begin{itemize}
\item[a)] $W = L$,
\item[b)] the natural map $K \rightarrow L$ is surjective,
\item[c)] the natural map $\chnum{X}/ h \cdot \chnum{X} \rightarrow 
\ol{\chow{*}(R)\subq}$ is an isomorphism.
\end{itemize}

In particular, we have the following theorem:

\bangou 
\begin{Theorem}\label{coin}
If ${\rm CH}^\cdot(X)\subq$ is isomorphic to 
$\chnum{X}$, then the natural map
$\chow{*}(R)\subq \rightarrow \ol{\chow{*}(R)\subq}$ is an isomorphism.
\end{Theorem}

In Roberts-Srinivas~\cite{RS}, the following statements are proved:
\begin{itemize}
\item[1)]
Suppose $k = {\Bbb C}$.
Then, there is an example such that $W \neq L$.
\item[2)]
Suppose $k = \ol{\Bbb Q}$ or $\ol{{\Bbb F}_p}$.
If some famous conjectures are true, then $W = L$.
\end{itemize}

Here, we give some examples.

\bangou
\begin{Example}\label{6.5.5}
\begin{rm}
Suppose that $n \geq 2$.
Let $X$ be the blowing-up 
of the projective space ${\Bbb P}^n_{\Bbb C}$
at $r$ distinct points.
Let $R$ be the local ring at the homogeneous maximal ideal 
of an affine cone of $X$.
Using Proposition~6.7 in Fulton~\cite{F}, it is easily verified that
\[
\chow{i}(X)\subq = \left\{
\begin{array}{lll}
{\Bbb Q} & & (i = 0,\ n) \\
{\Bbb Q}^{r+1} & & (i = 1, 2, \ldots, n-1) \\
0 & & (\mbox{otherwise})
\end{array}
\right.
\]
and ${\rm CH}^\cdot(X)\subq$ is isomorphic to 
$\chnum{X}$.
Then, the natural map
$\chow{*}(R)\subq \rightarrow \ol{\chow{*}(R)\subq}$ is an isomorphism
by Theorem~\ref{coin}.
Therefore, we obtain the following:
\[
\ol{\chow{i}(R)\subq} = \left\{
\begin{array}{lll}
{\Bbb Q} & & (i = n+1) \\
{\Bbb Q}^r & & (i = n) \\
0 & & (\mbox{otherwise})
\end{array}
\right.
\]
In this case, $\dim \ol{\g(R)\subq} = r+1$.
\end{rm}
\end{Example}

\bangou
\begin{Example}\label{6.6}
\begin{rm}
Let $m$ and $n$ be positive integers such that $n \geq m \geq 2$.
Suppose
\[
A = \left. k[x_{ij} \mid \mbox{$i = 1, \ldots, m$; $j = 1, \ldots, n$}]
\right/
I_2(x_{ij}) ,
\]
where $I_2(x_{ij})$ is the ideal of $A$ 
generated by all of the $2 \times 2$ minors of the $m \times n$ 
matrix $(x_{ij})$.
Then, $X = \proj(A) = {\Bbb P}^{m-1} \times {\Bbb P}^{n-1}$ and
$\dim R = m + n -1$.
In this case, ${\rm CH}^\cdot(X)\subq$ is isomorphic to 
$\chnum{X}$.
By Theorem~\ref{coin}, we have 
\[
\chow{i}(R)\subq = \ol{\chow{i}(R)\subq} =
\left\{
\begin{array}{cl}
{\Bbb Q} & (\mbox{if $i = n, n+1, \ldots, m+n-1$}) \\
0 & (\mbox{otherwise}) 
\end{array}
\right.
\]
and $\dim \ol{\g(R)\subq} = m$.
\end{rm}
\end{Example}

\bangou
\begin{Example}\label{6.7}
\begin{rm}
Set
\begin{eqnarray*}
A_{2n-1} & = & k[x_1,y_1, \ldots,x_n,y_n]/(x_1y_1 + \cdots + x_ny_n) \ \ 
\mbox{and} \\ 
A_{2n} & = & 
k[x_1,y_1, \ldots,x_n,y_n,z]/(z^2+x_1y_1 + \cdots + x_ny_n) .
\end{eqnarray*}
By Swan~\cite{Sw} (see Section~4 in Levine~\cite{Le}), we have
\[
{\rm CH}^i(\proj(A_{2n-1})) =
\left\{
\begin{array}{cll}
{\Bbb Q} & & (i = 0, 1, \ldots, n-2, n, \ldots, 2n-2) \\
{\Bbb Q}\oplus {\Bbb Q} & & (i = n-1) \\
0 & & (\mbox{otherwise})
\end{array}
\right.
\]
and
\[
{\rm CH}^i(\proj(A_{2n})) =
\left\{
\begin{array}{cll}
{\Bbb Q} & & (i = 0, 1, \ldots, 2n-1) \\
0 & & (\mbox{otherwise}) .
\end{array}
\right.
\]
Let $R_l$ denote the localization of $A_l$ at
the homogeneous maximal ideal.
Then, we can show
\[
\chow{i}(R_{2n-1})\subq = \ol{\chow{i}(R_{2n-1})\subq} =
\left\{
\begin{array}{cl}
{\Bbb Q} & (\mbox{if $i = n, 2n-1$}) \\
0 & (\mbox{otherwise}) 
\end{array}
\right.
\]
and
\[
\chow{i}(R_{2n})\subq = \ol{\chow{i}(R_{2n})\subq} =
\left\{
\begin{array}{cl}
{\Bbb Q} & (\mbox{if $i = 2n$}) \\
0 & (\mbox{otherwise}) 
\end{array}
\right.
\]
since ${\rm CH}^\cdot(\proj(A_l))\subq$ is isomorphic to 
$\chnum{\proj(A_l)}$ in this case.
\end{rm}
\end{Example}

\bangou
\begin{Example}\label{6.8}
\begin{rm}
Let $k$ be a field and set
\[
A = \left. k[x_{ij} \mid \mbox{$i = 1, \ldots, m$; $j = 1, \ldots, n$}]_{
(x_{ij} \mid i, \ j)}
\right/
I_t(x_{ij}) ,
\]
where $n \geq m \geq t$ are positive integers, and
$I_t(x_{ij})$ is the ideal of $A$ 
generated by all the $t \times t$ minors of the $m \times n$ matrix $(x_{ij})$.
Let $R$ be the local ring at the homogeneous maximal ideal of $A$.
Then, by Example~6.2 in \cite{K16}, the following three conditions are 
equivalent:
\begin{enumerate}
\item
$R$ is a Roberts ring,
\item
$R$ is a complete intersection, 
\item
$t = 1$ or $m = n = t$.
\end{enumerate}
In the case of $t = 2$, $R$ is a numerically Roberts ring if and only if
$m = n = 2$ (see Example~\ref{6.6}).
In the case of $t \geq 3$,
the author does not know when $R$ is a numerically Roberts ring.

Let $A_d(n)$ be the local ring at the homogeneous maximal ideal 
of the affine cone of the Grassmann variety $G_d(n)$ under the
Pl\"ucker embedding, where $d$ and $n$ are integers such that $0 < d < n$.
In \cite{K20}, it is proved that $A_d(n)$ is a Roberts ring 
if and only if one of the following conditions are satisfied,
\begin{enumerate}
\item $d=1$;
\item $d=n-1$;
\item $d=2$ and $n=4$;
\item $d=3$ and $n=6$.
\end{enumerate}
It is known that $\chring{G_d(n)} = \chnum{G_d(n)}$ 
(see Chapter~14 in \cite{F}).
Therefore, $A_d(n)$ is a numerically Roberts ring if and only if 
$A_d(n)$ is a Roberts ring.

Set 
\[
B_m(n) = \left. k[z_{ij} \mid 1 \leq i < j \leq n]_{
(z_{ij} \mid i, \ j)}
\right/
{\rm Pf}_m(z_{ij}) ,
\]
where ${\rm Pf}_m(z_{ij})$ is the pfaffian ideal of degree $m$ of
the $n \times n$ generic anti-symmetric matrix $(z_{ij})$
(see Section~5 in \cite{K20}).
By Theorem~5.1 in \cite{K20}, the following conditions are equivalent:
\begin{enumerate}
\item
$B_m(n)$ is a Roberts ring,
\item
$B_m(n)$ is a complete intersection, 
\item
$m = 1$ or $n = 2m $.
\end{enumerate}
In the case of $m = 2$, 
$B_2(n)$ coincides with $A_2(n)$.
Therefore, $B_2(n)$ is a numerically Roberts ring if and only if
$n = 4$.
In the case of $m \geq 3$,
the author does not know when $B_m(n)$ is a numerically Roberts ring.
\end{rm}
\end{Example}

\bangou
\begin{Remark}\label{6.9}
\begin{rm}
Let $R$ be the local ring at the homogeneous maximal ideal of
an affine cone $A$ of a smooth projective variety $X$ over ${\Bbb C}$.
Set $k = \dim X$.
Assume that the natural map 
\bangou
\begin{equation}\label{6.10}
\chhom{X} \longrightarrow \chnum{X}
\end{equation}
is an isomorphism,
where $\chhom{X}$ is the Chow ring (with rational coefficients) 
modulo homological equivalence
(Chapter~19 in \cite{F}).
Then, the fact that $\ol{\chow{j}(R)\subq} = 0$ for 
$j \leq (k + 1)/2 = \dim R/2$ can be proven
as follows.
(The map (\ref{6.10}) is an isomorphism in the case 
in which $\dim X$ is at most
$3$ or $X$ is an abelian variety (Example~19.3.2 in \cite{F}).
If we assume that Grothendieck's standard conjectures are true,
then the map (\ref{6.10}) is an isomorphism.)

We have an injective map
\[
cl^X : {\rm CH}_{\rm hom}^i(X)\subq \longrightarrow H^{2i}(X,{\Bbb Q})
\]
for $i = 0, 1, \ldots, k$ called the {\em cycle map} (see Chapter~19 in
\cite{F}).
Let $h$ be the very ample divisor of the embedding $X = \proj(A)$.
Since the map (\ref{6.10}) is an isomorphism,
we have the following commutative diagram for $i > k/2$:
\[
\begin{array}{ccc}
{\rm CH}_{\rm num}^{k-i}(X)\subq & \stackrel{cl^X}{\ya{aaaaaa}} & 
H^{2(k-i)}(X,{\Bbb Q})  \\
{\scriptstyle h^{2i-k}} \downarrow 
\hphantom{\scriptstyle h^{2i-k}} & & 
{\scriptstyle h^{2i-k}} \downarrow 
\hphantom{\scriptstyle h^{2i-k}} \\
{\rm CH}_{\rm num}^i(X)\subq & \stackrel{cl^X}{\ya{aaaaaa}}
 & H^{2i}(X,{\Bbb Q})
\end{array}
\]
Note that the horizontal maps are injective.
By the hard Lefschetz theorem (cf.\ \cite{Le}),
the right vertical map in the above diagram is an isomorphism.
Since $\dim {\rm CH}_{\rm num}^{k-i}(X)\subq = \dim
{\rm CH}_{\rm num}^i(X)\subq$,
the left vertical map is also an isomorphism.
In particular, the map
\[
{\rm CH}^{i-1}_{\rm num}(X)\subq \stackrel{h}{\longrightarrow} 
{\rm CH}^i_{\rm num}(X)\subq 
\]
is surjective for $i > k/2$.

Recall that there exist surjective maps
\[
{\rm CH}_{\rm num}^i(X)\subq / h \cdot {\rm CH}_{\rm num}^{i-1}(X)\subq
\rightarrow 
\ol{\chow{\dim R - i}(R)\subq}
\]
for each $i$ by (\ref{6.4}).
Therefore, we have $\ol{\chow{j}(R)\subq} = 0$ for $j \leq (k+1)/2 = \dim R/2$.
By Remark~\ref{3.4} and (\ref{7.7}), this is equivalent to the condition that
$\chi_{{\Bbb F}.}(M)$ is equal to $0$ for any ${\Bbb F}. \in C^{\m}(R)$
and any finitely generated $R$-module $M$ with $\dim M \leq \dim R/2$.
\end{rm}
\end{Remark}

\section{A vanishing of Chow group modulo numerical equivalence}\label{7}
In this section, we attempt to prove the following theorem:

\bangou
\begin{Theorem}\label{7.1}
\vanishing
\end{Theorem}

By this theorem, we know that, for any regular alteration
$\pi' : Z' \rightarrow \spec(R)$, we have 
\[
\dim \pi'^{-1}(\m) \geq d - \min\{ t \mid \ol{\chow{t}(R)\subq} \neq 0 \} .
\]
Before proving the theorem, we give some examples.

\bangou
\begin{Example}\label{7.2}
\begin{rm}
Let $(R,\m)$ be a $d$-dimensional Noetherian local domain
and let $f : W \rightarrow \spec(R)$ be a resolution of singularities
of $\spec(R)$.
Then, by the above argument, we have 
\[
\dim f^{-1}(\m) \geq d-\min\{ t \mid \ol{\chow{t}(R)\subq} \neq 0 \} .
\]

Let $R$ be a ring in Example~\ref{6.6}.
In this case, $d$ is equal to $m + n -1$ and we have 
\[
\min\{ t \mid \ol{\chow{t}(R)\subq} \neq 0 \} = n .
\]
On the other hand, let $g : V \rightarrow \spec(R)$ be the blowing-up
of $\spec(R)$ with center $(x_{11}, x_{21}, \ldots, x_{m1})$.
Then, $g$ is a resolution of singularities of $\spec(R)$ with
$\dim g^{-1}(\m) = m-1$.
Therefore, in this case, 
\[
\dim g^{-1}(\m) = d-\min\{ t \mid \ol{\chow{t}(R)\subq} \neq 0 \} .
\]

Next, we assume that $R$ is the local ring at the homogeneous maximal ideal
of an affine cone of an abelian variety $X$ over the field of complex numbers.
Then, it is known that, for any regular alteration 
$\pi : Z \rightarrow \spec(R)$, $\dim \pi^{-1}(\m)$ is equal to 
$\dim X = d-1$ (cf.\ Example~3.9 in Ishii-Milman~\cite{IM}).
On the other hand, by Remark~\ref{6.9}, we have
\[
\min\{ t \mid \ol{\chow{t}(R)\subq} \neq 0 \} > d/2 .
\]
Therefore, if $d \geq 3$, we have
\[
\dim \pi^{-1}(\m) > d-\min\{ t \mid \ol{\chow{t}(R)\subq} \neq 0 \}
\]
for any regular alteration $\pi : Z \rightarrow \spec(R)$.
\end{rm}
\end{Example}

In order to prove Theorem~\ref{7.1}, we need the following lemma:

\bangou
\begin{Lemma}\label{7.3}
Let $X$ be a $d$-dimensional integral scheme that is of finite type over a
regular scheme.
Assume that there is a regular alteration 
$\pi : Z \rightarrow X$.
Let $Y$ be a closed subset of $X$. 
Let $i$ be a non-negative integer and take
$c \in {\rm A}^i(Y \rightarrow X)\subq$.
Then, the following conditions are equivalent:
\begin{itemize}
\item[1)]
$c = 0$ in ${\rm A}^i(Y \rightarrow X)\subq$.
\item[2)]
$c([Z]) = 0$ in $\chow{d-i}(\pi^{-1}(Y))\subq$.
\end{itemize}
In particular, if $d-i > \dim \pi^{-1}(Y)$,
then ${\rm A}^i(Y \rightarrow X)\subq$ coincides with $0$.
\end{Lemma}

Before proving Lemma~\ref{7.3}, recall bivariant classes and 
the bivariant group ${\rm A}^i(Y \rightarrow X)$
defined in Chapter~17 in Fulton~\cite{F}.
A bivariant class $c$ in ${\rm A}^i(Y \rightarrow X)$ is a collection 
of homomorphisms
\[
c^{(k)}_g : \chow{k}(X') \longrightarrow \chow{k-i}(X' \times_X Y)
\]
for all $g : X' \rightarrow X$ and all $k \in {\Bbb Z}$
compatible with proper push-forward, flat pull-back, and intersection products,
that is, $c^{(k)}_g$'s satisfy (${\rm C}_1$), (${\rm C}_2$) and 
(${\rm C}_3$) in Definition~17.1 in \cite{F}.

\vspace{2mm}

\proof
By definition, 1) implies 2).

We shall prove that 2) implies 1).
Let $V$ be an $s$-dimensional integral scheme
and let $g : V \rightarrow X$ be a morphism
of finite type.
We want to show $c([V]) = 0$ in $\chow{s-i}(g^{-1}(Y))\subq$.

Take a closed integral subscheme $V'$ of $V \times_XZ$ such that
the composite map $V' \rightarrow V \times_XZ \rightarrow V$
is proper surjective and generically finite.
Then, by the compatibility with proper push-forward 
(Definition~17.1 (${\rm C}_1$) in \cite{F}),
$c([V']) = 0$ implies that $c([V]) = 0$.
Replacing $V$ by $V'$, we may assume that there is a morphism
$h : V \rightarrow Z$ such that $\pi \cdot h = g$.

Setting $W = \pi^{-1}(Y)$,
we have a natural map
\[
{\rm A}^i(Y \rightarrow X)\subq \rightarrow {\rm A}^i(W \rightarrow Z)\subq 
\]
as in 17.2 (${\rm P}_3$) in \cite{F}.
The image of $c$ will be denoted by $c$ again.
Then, $c \in {\rm A}^i(W \rightarrow Z)\subq$ satisfies $c([Z]) = 0$
by the assumption.
In this situation,
we want to show that $c([V]) = 0$.

By Nagata's compactification~(\cite{Lu}, \cite{N}),
there exist an integral scheme $\ol{V}$,
a proper morphism $\ol{h} : \ol{V} \rightarrow Z$ and
an open immersion $j : V \rightarrow \ol{V}$
such that $h = \ol{h} \cdot j$.
Then, by the compatibility with flat pull-back 
(Definition~17.1 (${\rm C}_2$) in \cite{F}),
$c([\ol{V}]) = 0$ implies $c([V]) = 0$.
Replacing $V$ by $\ol{V}$, we may assume that the morphism
$h : V \rightarrow Z$ is proper.

Furthermore, using Chow's lemma (cf.\ \cite{Ha}),
we may assume that $h : V \rightarrow Z$ is projective 
by the compatibility with proper push-forward.
Suppose that $V$ is a closed subscheme of ${\Bbb P}^n_Z$
and the composite map 
\[
V \subseteq {\Bbb P}^n_Z \rightarrow Z
\]
coincides with $h$.
Set $t = \codim_{{\Bbb P}^n_Z}V$.
Since ${\Bbb P}^n_Z$ is a regular scheme,
there exists a locally-free ${\cal O}_{{\Bbb P}^n_Z}$-resolution
${\Bbb F}.$ of ${\cal O}_V$.
Then, we have 
\[
{\rm ch}_t({\Bbb F}.)([{\Bbb P}^n_Z]) = [V]
\]
by the top term property~(Theorem~18.3~(5) in \cite{F}).
Therefore, we have
\begin{eqnarray*}
c([V]) & = & c({\rm ch}_t({\Bbb F}.)([{\Bbb P}^n_Z])) \\
& = & {\rm ch}_t({\Bbb F}.)(c([{\Bbb P}^n_Z])) 
\end{eqnarray*}
since $c$ and ${\rm ch}_t({\Bbb F}.)$ commute by a theorem of 
Roberts~\cite{R1}.
Since ${\Bbb P}^n_Z \rightarrow Z$ is a flat map,
$c([Z]) = 0$ implies $c([{\Bbb P}^n_Z]) = 0$ by the compatibility 
with flat pull-back.
Therefore, we have $c([V]) = 0$.
\qed

\bangou
\begin{Corollary}\label{7.4}
Let $X$ be a $d$-dimensional integral scheme that is
of finite type over a regular scheme.
Assume that there exists a regular alteration 
$\pi : Z \rightarrow X$.
Let $Y$ be a closed subset of $X$ and set $W = \pi^{-1}(Y)$.
Let ${\Bbb F}.$ be a bounded ${\cal O}_X$-locally free complex
with support in $Y$.
Then, for a non-negative integer $k$, the following conditions are equivalent:
\begin{itemize}
\item[1)]
${\rm ch}_i({\Bbb F}.) = 0$ in ${\rm A}^i(Y \rightarrow X)\subq$
for $i = 0, 1, \ldots, k-1$,
\item[2)]
$\sum_i(-1)^i[H_i(\pi^*{\Bbb F}.)] \in F_{d-k}\g(W)\subq$,
\end{itemize}
where $F_{d-k}\g(W)\subq$ is a vector subspace of $\g(W)\subq$
spanned by 
\[
\{ [{\cal M}] \mid \mbox{${\cal M}$ is a coherent
${\cal O}_W$-module with $\dim \supp {\cal M} \leq d-k$} \} .
\]
\end{Corollary}

Note that condition 1) above is independent of
the choice of a regular alteration $\pi : Z \rightarrow X$.

\vspace{2mm}

\proof
By Lemma~\ref{7.3}, condition 1) is equivalent to the condition that
${\rm ch}_i(\pi^*{\Bbb F}.)([Z]) = 0$ in $\chow{d-i}(W)\subq$
for $i = 0, 1, \ldots, k-1$.
Furthermore, this is equivalent to 
\[
{\rm ch}(\pi^*{\Bbb F}.)([Z])
\in \bigoplus_{j = 0}^{d-k} \chow{j}(W)\subq .
\]

Since the diagram
\[
\begin{array}{ccc}
\g(Z)\subq & \stackrel{\tau_{Z/Z}}{\longrightarrow} & \chow{*}(Z)\subq \\
{\scriptstyle \chi_{\pi^*{\Bbb F}.}} \downarrow
\hphantom{\scriptstyle \chi_{\pi^*{\Bbb F}.}} & & 
{\scriptstyle {\rm ch}(\pi^*{\Bbb F}.)} \downarrow
\hphantom{\scriptstyle {\rm ch}(\pi^*{\Bbb F}.)} \\
\g(W)\subq & \stackrel{\tau_{W/Z}}{\longrightarrow} & \chow{*}(W)\subq 
\end{array}
\]
is commutative by the local Riemann-Roch theorem~\cite{F},
we have 
\[
\tau_{W/Z}\left( \sum_i(-1)^i[H_i(\pi^*{\Bbb F}.)] \right) =
{\rm ch}(\pi^*{\Bbb F}.)(\tau_{Z/Z}([{\cal O}_Z])) =
{\rm ch}(\pi^*{\Bbb F}.)([Z]) .
\]
Recall that, since $Z$ is a regular base scheme, we have 
$\tau_{Z/Z}([{\cal O}_Z]) = [Z]$.

By the top term property, we have 
$\tau_{W/Z}(F_{d-k}\g(W)\subq) = \oplus_{j = 0}^{d-k} \chow{j}(W)\subq$.
Therefore, 1) is equivalent to 2).
\qed

Next, we define two invariants for complexes in $C^\m(R)$.

\bangou
\begin{Definition}\label{7.5}
\begin{rm}
Let $(R,\m)$ be a Noetherian local ring of dimension $d$.
For a complex ${\Bbb F}. \in C^\m(R)$, we define
\begin{eqnarray*}
b({\Bbb F}.) & = & \min\{
s \mid \mbox{${\rm ch}_s({\Bbb F}.) \neq 0$ in 
${\rm A}^s(\spec(R/\m) \rightarrow \spec(R))\subq$} \} , \\
n({\Bbb F}.) & = & \min\{
s \mid \mbox{${\rm ch}_s({\Bbb F}.) : \chow{s}(R)\subq \rightarrow {\Bbb Q}$ 
does not coincide with $0$} \} .
\end{eqnarray*}
If ${\rm ch}_s({\Bbb F}.)$ is $0$ in 
${\rm A}^s(\spec(R/\m) \rightarrow \spec(R))\subq$ 
for all $s$,
we set $b({\Bbb F}.) = \infty$.
If the map ${\rm ch}_s({\Bbb F}.) : \chow{s}(R)\subq \rightarrow {\Bbb Q}$ 
is equal $0$ for all $s$, 
we set $n({\Bbb F}.) = \infty$.
\end{rm}
\end{Definition}

By definition, either $0 \leq n({\Bbb F}.) \leq d$ or $n({\Bbb F}.) = \infty$.

It is easily verified that
\bangou
\begin{equation}\label{7.6}
0 \leq b({\Bbb F}.) \leq n({\Bbb F}.) .
\end{equation}
Furthermore, note that
\bangou
\begin{equation}\label{7.7}
\min\{ t \mid \ol{\chow{t}(R)\subq} \neq 0 \}
= \min\{ n({\Bbb F}.) \mid {\Bbb F}. \in C^\m(R) \} .
\end{equation}

\bangou
\begin{Remark}\label{7.8}
\begin{rm}
Let $(R,\m)$ be a $d$-dimensional 
Noetherian local domain with regular alteration
$\pi : Z \rightarrow \spec(R)$.
Set $X = \pi^{-1}(\m)$ and 
$\chi(\pi^*{\Bbb F}.) = \sum_i(-1)^i[H_i(\pi^*{\Bbb F}.)]
\in \g(X)\subq$.
Then, for each ${\Bbb F}. \in C^\m(R)$, we have 
\[
\chi(\pi^*{\Bbb F}.) \in F_{d - b({\Bbb F}.)}\g(X)\subq \setminus
F_{d - b({\Bbb F}.) -1}\g(X)\subq 
\]
by Corollary~\ref{7.4}.
This is equivalent to
\bangou
\begin{equation}\label{7.9}
\tau_{X/Z}(\chi(\pi^*{\Bbb F}.)) \in
\bigoplus_{i \leq d - b({\Bbb F}.)}\chow{i}(X)\subq
\setminus \bigoplus_{i < d- b({\Bbb F}.)}\chow{i}(X)\subq .
\end{equation}
In particular,
either $0 \leq b({\Bbb F}.) \leq d$ or $b({\Bbb F}.) = \infty$.
\end{rm}
\end{Remark}

Now we return to the proof of Theorem~\ref{7.1}.
Let $(R, \m)$ be a local ring that satisfies the assumptions in
Theorem~\ref{7.1}.
Then, for ${\Bbb F}. \in C^\m(R)$, 
we have
\[
\sum_i(-1)^i[H_i(\pi^*{\Bbb F}.)] \in 
F_{d-b({\Bbb F}.)}\g(\pi^{-1}(\m))\subq \setminus
F_{d-b({\Bbb F}.)-1}\g(\pi^{-1}(\m))\subq .
\]
In particular, we have
\[
d-b({\Bbb F}.) \leq \dim \pi^{-1}(\m) 
\]
for any ${\Bbb F}. \in C^\m(R)$.
It follows that 
\[
d-\dim \pi^{-1}(\m) \leq 
\min\{ b({\Bbb F}.) \mid {\Bbb F}. \in C^\m(R) \} \leq
\min\{ n({\Bbb F}.) \mid {\Bbb F}. \in C^\m(R) \} =
\min\{ t \mid \ol{\chow{t}(R)\subq} \neq 0 \}
\]
by (\ref{7.6}) and (\ref{7.7}).

We have completed the proof of Theorem~\ref{7.1}.

\bangou
\begin{Remark}\label{7.10}
\begin{rm}
Let $(R, \m)$ be a $d$-dimensional Noetherian local ring.

Let ${\Bbb K}.$ be a Koszul complex with respect to a system of parameters
for $R$.
Then, $b({\Bbb K}.) = n({\Bbb K}.) = d$.

Suppose that $R$ is a homomorphic image of a regular local ring $(S,\n)$.
Assume that ${\Bbb F}. \in C^\m(R)$ is liftable to $S$, that is,
there exists an $S$-free complex ${\Bbb G}.$ such that
${\Bbb F}. = {\Bbb G}.\otimes_SR$.
(Note that a Koszul complex of a system of parameters for $R$
is always liftable.) \
Then,
\[
b({\Bbb F}.) = \mbox{$d$ or $\infty$} .
\]
In fact, let $Y$ be the support of the complex ${\Bbb G}.$, i.e.,
$Y = \cup_i\supp(H_i({\Bbb G}.)) \subset \spec(S)$.
Then, $\spec(R) \cap Y = \{ \n \}$.
Since $\dim R + \dim Y \leq \dim S$ by Serre's theorem~\cite{Se},
we have $\dim Y \leq \dim S - d$.
Since $\spec(S)$ itself is a regular alteration of $\spec(S)$,
we have ${\rm ch}_i({\Bbb G}.) = 0$
in ${\rm A}^i(Y \rightarrow \spec(S))\subq$
for $i = 0, 1, \cdots, d-1$ by Corollary~\ref{7.4}.
It follows that ${\rm ch}_i({\Bbb F}.)$ coincides with $0$
in ${\rm A}^i(\spec(R/\m) \rightarrow \spec(R))\subq$
for $i = 0, 1, \cdots, d-1$.
\end{rm}
\end{Remark}

\bangou
\begin{Remark}\label{7.11}
\begin{rm}
We discuss affine cones of smooth projective varieties.
With notation as in Section~\ref{6},
we have the following commutative diagram:
\[
\begin{array}{ccccc}
\k{\m}(R)\subq & \stackrel{\varphi}{\longrightarrow} & K
& \subseteq & {\rm CH}^\cdot(X)\subq \\
& {\scriptstyle \phi} \searrow \hphantom{\scriptstyle \phi} & 
\downarrow & & \downarrow \\
& & W & \subseteq & \chnum{X}
\end{array}
\]
Then, by the definition of $\varphi$ (see (\ref{6.2})) and Remark~\ref{7.8},
we have
\[
\varphi([{\Bbb F}.])
\in \bigoplus_{i \geq b({\Bbb F}.)-1}{\rm CH}^i(X)\subq
\setminus \bigoplus_{i > b({\Bbb F}.)-1}{\rm CH}^i(X)\subq
\]
for each ${\Bbb F}. \in C^\m(R)$.

On the other hand, using the method described 
in Section~\ref{6}, we can prove
\[
\phi([{\Bbb F}.])
\in \bigoplus_{i \geq n({\Bbb F}.)-1}{\rm CH}^i_{\rm num}(X)\subq
\setminus \bigoplus_{i > n({\Bbb F}.)-1}{\rm CH}^i_{\rm num}(X)\subq .
\]
Therefore, the difference between $b({\Bbb F}.)$ and $n({\Bbb F}.)$
corresponds to the difference between rational equivalence and
numerical equivalence on $X$.
\end{rm}
\end{Remark}

\bangou
\begin{Example}\label{7.11.5}
\begin{rm}
Let $R$ be a ring as in Example~\ref{6.6}.
In this case, since ${\rm CH}^\cdot(X)\subq$ is isomorphic to $\chnum{X}$,
$b({\Bbb F}.)$ coincides with $n({\Bbb F}.)$ for any ${\Bbb F}. \in C^\m(R)$.
By (\ref{7.7}), 
we have $b({\Bbb F}.) = n({\Bbb F}.) \geq n$.

Let ${\Bbb H}\{ s \}.$ be the complex as in Remark~\ref{2.4}.
Then, $b({\Bbb H}\{ s \}.) = n({\Bbb H}\{ s \}.) = s$
for $s = n, n+1, \ldots, m+n-1$.
\end{rm}
\end{Example}

\bangou
\begin{Example}\label{7.12}
\begin{rm}
We give an example of ${\Bbb F}. \in C^\m(R)$
with $b({\Bbb F}.) \neq n({\Bbb F}.)$.

Set
\[
A = {\Bbb C}[x, y, z]/(x^3+y^3+z^3) .
\]
We regard $A$ as a graded ring with $\deg(x) = \deg(y) = \deg(z) = 1$.
Set $X = \proj(A)$,
and let $(R,\m)$ be the local ring of $A$ at $(x,y,z)$.
With notation as in Section~\ref{6},
the map 
\[
\k{\m}(R)\subq \stackrel{\varphi}{\longrightarrow} 
K = {\rm Ker}({\rm CH}^\cdot(X)\subq \stackrel{h}{\longrightarrow}
{\rm CH}^\cdot(X)\subq)
\]
is surjective.
Take $0 \neq \alpha \in {\rm CH}^1(X)\subq$ such that 
$0 = \ol{\alpha} \in {\rm CH}^1_{\rm num}(X)\subq$.
(Since $X$ is an elliptic curve over ${\Bbb C}$,
it is easy to take such $\alpha$.)
Then, by the surjectivity of $\varphi$,
there exists ${\Bbb F}. \in C^\m(R)$ such that
$\varphi([{\Bbb F}.]) = n\alpha$ with some positive integer $n$.
Then, by Remark~\ref{7.11},
we have 
\[
b({\Bbb F}.) = 2 \neq \infty = n({\Bbb F}.)  .
\]
\end{rm}
\end{Example}

\section{Applications of Lemma~\ref{7.3}}\label{8}

In this section, we give two applications of Lemma~\ref{7.3}.
The first is a sufficient condition of the vanishing property 
of intersection 
multiplicities~(see Theorem~\ref{8.1} below). 
The second gives another proof of the vanishing theorem 
(see Theorem~\ref{8.5}
below) of the first localized Chern characters due to Roberts~\cite{Rmac}.

\bangou
\begin{Theorem}\label{8.1}
\intersectionmultiplicity
\end{Theorem}

\proof
Let $M$ and $N$ be $R$-modules that satisfy the assumptions in the theorem.
We may assume $\dim M \geq \dim N$.
Then, $\dim N < d/2$.
Let ${\Bbb F}.$ and ${\Bbb G}.$ be finite free resolutions of
$M$ and $N$, respectively.

Set $\tau_{R/S}([R]) = \tau_d + \tau_{d-1} + \cdots + \tau_0$, where
$\tau_i \in \chow{i}(R)\subq$.
Then, by the local Riemann-Roch formula, we have
\begin{eqnarray*}
\sum_i(-1)^i\ell_R(\tor{i}{R}{M}{N}) & = & 
{\rm ch}({\Bbb F}.\otimes{\Bbb G}.) \cap \tau_{R/S}([R]) \\
& = & \sum_{i,\ j \geq 0} {\rm ch}_i({\Bbb F}.){\rm ch}_j({\Bbb G}.) \cap 
\tau_{i + j} .
\end{eqnarray*}
We want to prove ${\rm ch}_i({\Bbb F}.){\rm ch}_j({\Bbb G}.) \cap 
\tau_{i + j} = 0$ for any $i, j \geq 0$.

Assume that ${\rm ch}_i({\Bbb F}.){\rm ch}_j({\Bbb G}.) \cap 
\tau_{i + j} \neq 0$ for some $i, j \geq 0$.
By Lemma~\ref{7.3}, we have
\bangou
\begin{equation}\label{8.2}
\dim \pi^{-1}(\supp(M)) \geq d-i
\end{equation}
since ${\rm ch}_i({\Bbb F}.) \neq 0$ in 
${\rm A}^i(\supp(M) \rightarrow \spec(R))\subq$.
%
Furthermore, since
\[
0 \neq {\rm ch}_j({\Bbb G}.) \cap \tau_{i + j} \in \chow{i}(\supp(N))\subq ,
\]
we have $0 \leq i \leq \dim N$ and 
\[
\dim M + i \leq \dim M + \dim N < d .
\]
Therefore,
\[
\dim M < d-i .
\]

We have
\begin{eqnarray*}
\pi^{-1}(\supp(M)) \subseteq \pi^{-1}(Y \cup \supp(M))
& = & \pi^{-1}(Y) \cup \pi^{-1}(\supp(M) \setminus Y) \\
& = & \pi^{-1}(Y) \cup \ol{\pi^{-1}(\supp(M) \setminus Y)} ,
\end{eqnarray*}
where $\ol{\pi^{-1}(\supp(M) \setminus Y)}$ denotes the closure
of $\pi^{-1}(\supp(M) \setminus Y)$.
It is easy to see that
\[
\dim \ol{\pi^{-1}(\supp(M) \setminus Y)} \leq \dim M ,
\]
since $R$ is an excellent ring.
Therefore, we have
\[
\dim \pi^{-1}(\supp(M)) \leq \max\{ \dim \pi^{-1}(Y), \dim M \} .
\]

Since $i \leq \dim N < d/2$, 
we obtain
\[
\dim \pi^{-1}(Y) \leq d/2 < d-i .
\]
Then, we have 
\[
\dim \pi^{-1}(\supp(M)) \leq \max\{ \dim \pi^{-1}(Y), \dim M \}
< d-i .
\]
This inequality contradicts (\ref{8.2}).
\qed

\bangou
\begin{Example}
\begin{rm}
Let $m$, $n$, $q$ be positive integers such that $2 \leq m \leq n$.
Set
\[
R = \left.
k\left[
\{ x_{ij} \mid i = 1, \cdots, m; \ j = 1, \ldots, n \}, \ 
y_1, \ldots, y_q
\right]_{(\underline{x}, \underline{y})}
\right/
I_2(x_{ij}) .
\]
Then, we have $\dim R = m + n + q -1$.
Let $\pi : Z \rightarrow \spec(R)$ be the blowing-up of $\spec(R)$
along $(x_{11}, x_{21}, \ldots, x_{m1})$.
Then, $Z$ is a resolution of singularities of $\spec(R)$
such that $\dim \pi^{-1}((x_{11}, x_{21}, \ldots, x_{m1})) = m+q-1$.
Therefore, if $m+q \leq n+1$, then $\pi : Z \rightarrow \spec(R)$ satisfies
the assumptions in Theorem~\ref{8.1}.

Note that the dimension of the singular locus of $\spec(R)$
is equal to $q$.
(For a local ring such that the dimension of its singular locus is at most
$1$, the vanishing property was proved by Roberts~\cite{R1}.)
\end{rm}
\end{Example}

In the remainder of this section, we shall give another proof to
the vanishing theorem of the first localized Chern characters
due to Roberts~\cite{Rmac} as follows:

\bangou
\begin{Theorem}\label{8.5}
Let $X$ be a scheme that is of finite type over an excellent regular scheme.
Let $X = X_1 \cup \cdots \cup X_t$ be the irreducible decomposition.
Assume that each $(X_i)_{\rm red}$ has a regular alteration.
Let $Y$ be a closed subset of $X$ such that 
${\rm codim}_{X_i}X_i \cap Y \geq 2$ for $i = 1, \ldots, t$.
Then, for each ${\Bbb F}. \in C^Y(X)$, we have 
\[
{\rm ch}_1({\Bbb F}.) = 0
\]
in ${\rm A}^1(Y \rightarrow X)\subq$.
\end{Theorem}

Note that Roberts proved the vanishing theorem of 
the first localized Chern 
characters without assuming that 
$(X_i)_{\rm red}$'s have regular alterations.

\vspace{2mm}

\proof
Let $V$ be an $s$-dimensional integral scheme with a morphism
$f : V \rightarrow X$ of finite type.
We want to show that ${\rm ch}_1(f^*{\Bbb F}.) : \chow{s}(V)\subq
\rightarrow \chow{s-1}(f^{-1}(Y))\subq$ coincides with $0$.

Since $V$ is irreducible, some irreducible component $X_i$ contains $f(V)$.
Let $f_i : V \rightarrow X_i$ be the induced map.
We have only to show that 
\[
{\rm ch}_1(f_i^*({\Bbb F}.\otimes_{{\cal O}_X}{\cal O}_{X_i})) : 
\chow{s}(V)\subq \rightarrow \chow{s-1}(f_i^{-1}(X_i \cap Y))\subq
\]
coincides with $0$.
Therefore, we may assume that $X$ is an integral scheme
with a regular alteration.

Set $d = \dim X$.
Let $\pi : Z \rightarrow X$ be a regular alteration.
By Proposition~18.1~(a) in \cite{F}, 
${\rm ch}_0({\Bbb F}.)$ is the multiplication by 
$\sum_i(-1)^i\rank_{{\cal O}_X} F_i$.
Since $Y$ is a proper closed subset of $X$, we have
${\rm ch}_0({\Bbb F}.) = 0$ in ${\rm A}^0(Y \rightarrow X)\subq$.
Set $W = \pi^{-1}(Y)$.
By Corollary~\ref{7.4},
${\rm ch}_1({\Bbb F}.) = 0$ in ${\rm A}^1(Y \rightarrow X)\subq$
if and only if
\[
\tau_{W/Z}(\chi_{\pi^*{\Bbb F}.}([{\cal O}_Z]))
\in \bigoplus_{i=0}^{d-2}\chow{i}(W)\subq .
\]

Since $W$ is a proper closed subset of $Z$, we have $\dim W \leq d-1$.
If $\dim W \leq d-2$, we have nothing to prove.

Assume that $\dim W = d-1$. 
Let $W_1$, \ldots, $W_l$ be the irreducible components 
of $W$ of dimension $d-1$.
Then, $\chow{d-1}(W)\subq$ is a ${\Bbb Q}$-vector space with basis
$[W_1]$, \ldots, $[W_l]$.
It is sufficient to show that the coefficient of $[W_i]$ in 
$\tau_{W/Z}(\chi_{{\pi^*}{\Bbb F}.}([{\cal O}_Z]))$ is equal to $0$
for $i = 1, \ldots, l$.

Assume that the coefficient of $[W_1]$ in 
$\tau_{W/Z}(\chi_{\pi^*{\Bbb F}.}([{\cal O}_Z]))$ is not $0$.

Let $U$ be an affine open set of $X$ such that
$F_i|_U$ is ${\cal O}_U$-free for each $i$,
and $\pi^{-1}(U) \cap W_1 \neq \emptyset$.
Consider the following commutative diagram:
\[
\begin{array}{ccc}
\k{Y}(X)\subq & \longrightarrow & \k{U \cap Y}(U)\subq \\
{\scriptstyle \pi^*} \downarrow \phantom{\scriptstyle \pi^*} & & 
{\scriptstyle {\pi_U}^*} \downarrow \phantom{\scriptstyle {\pi_U}^*} \\
\k{W}(Z)\subq & \longrightarrow & \k{\pi^{-1}(U) \cap W}(\pi^{-1}(U))\subq \\
{\scriptstyle \chi} \downarrow \phantom{\scriptstyle \chi} & & 
{\scriptstyle \chi} \downarrow \phantom{\scriptstyle \chi} \\
\g(W)\subq & \longrightarrow & \g(\pi^{-1}(U) \cap W)\subq \\
{\scriptstyle \tau_{W/Z}} \downarrow \phantom{\scriptstyle \tau_{W/Z}} & & 
{\scriptstyle \tau_{\pi^{-1}(U) \cap W/\pi^{-1}(U)}} \downarrow 
\phantom{\scriptstyle \tau_{\pi^{-1}(U) \cap W/\pi^{-1}(U)}} \\
\chow{*}(W)\subq & \longrightarrow & \chow{*}(\pi^{-1}(U) \cap W)\subq
\end{array}
\]
The horizontal maps are restrictions to $U$.
By the commutativity, the coefficient of $[\pi^{-1}(U) \cap W_1]$ in 
\[
\tau_{\pi^{-1}(U) \cap W/\pi^{-1}(U)}
\chi {\pi_U}^*  ([{\Bbb F}.|_U])
\]
is equal to the coefficient of $[W_1]$ in 
$\tau_{W/Z}(\chi_{\pi^*{\Bbb F}.}([{\cal O}_Z]))$.

Replacing $X$ by $U$, we may assume that $X$ is an affine integral scheme
and ${\Bbb F}.$ is a complex of free modules.

Setting $\ol{X} = \spec(H^0(Z, {\cal O}_Z))$, we take the Stein 
factorization
\[
\begin{array}{clc}
Z & & \\
\downarrow &  \searrow{\scriptstyle \pi} & \\
\ol{X} & \stackrel{g}{\longrightarrow} & X
\end{array}
\]
of $\pi : Z \rightarrow X$.

Since $X$ is an excellent scheme, $\codim_{\ol{X}}g^{-1}(Y)$ is at least $2$.
Therefore, replacing $X$ by $\ol{X}$, 
we may assume that $X$ is an affine normal scheme
and $\pi : Z \rightarrow X$ is birational.

We denote by $i$ the inclusion $W \rightarrow Z$.
Consider the following commutative diagram:
\[
\begin{array}{ccccccc}
\k{Y}(X)\subq & \stackrel{\pi^*}{\ya{aaaa}} & \k{W}(Z)\subq & 
\stackrel{\chi}{\ya{aaaa}} &
\g(W)\subq & \stackrel{i_*}{\ya{aaaa}} & \g(Z)\subq \\
& & & & {\scriptstyle \tau_{W/Z}}\downarrow\hphantom{\scriptstyle \tau_{W/Z}}
& & \hphantom{\scriptstyle \tau_{Z/Z}}\downarrow{\scriptstyle \tau_{Z/Z}} \\
& & & & \chow{*}(W)\subq & \stackrel{i_*}{\ya{aaaa}} & \chow{*}(Z)\subq \\
& & & & {\scriptstyle p}\downarrow\hphantom{\scriptstyle p} & & 
\hphantom{\scriptstyle p}\downarrow{\scriptstyle p} \\
& & & & \chow{d-1}(W)\subq & \stackrel{i_*}{\ya{aaaa}} & \chow{d-1}(Z)\subq
\end{array}
\]
where $i_*$ is induced by the proper morphism $i : W \rightarrow Z$
and $p$'s are the projections.
Set $\alpha_{d-1} = p \cdot \tau_{W/Z} \cdot \chi(\pi^*{\Bbb F}.)$.
By our assumption, we have $\alpha_{d-1} \neq 0$.

On the other hand, we have
\[
i_*\chi(\pi^*{\Bbb F}.) = \sum_j(-1)^j[\pi^*F_j] =
\left( \sum_j(-1)^j\rank F_j \right) \cdot [{\cal O}_Z]= 0 ,
\]
because each $F_j$ is a free ${\cal O}_X$-module and $Y$ is 
a proper closed subset of $X$.
Therefore, $i_*(\alpha_{d-1}) = 0$.
This equality contradicts the following claim:

\bangou
\begin{Claim}\label{8.6}
The map $i_* : \chow{d-1}(W)\subq \longrightarrow \chow{d-1}(Z)\subq$ is injective.
\end{Claim}

In the remainder of this section, we shall prove the claim.

Let $W_1$, \ldots, $W_l$ be the irreducible components of $W$ of dimension $d-1$.
Then, $\chow{d-1}(W)\subq$ is a ${\Bbb Q}$-vector space with basis
$[W_1]$, \ldots, $[W_l]$.
Suppose that $i_*(q_1[W_1] + \cdots + q_l[W_l]) = 0$.
Then, there exists a positive integer $n$ and $a \in R(Z)^{\mathop{\times}} =
R(X)^{\mathop{\times}}$ such that
\[
{\rm div}_Z(a) = nq_1[W_1] + \cdots + nq_l[W_l] .
\]
Then, we have 
\[
{\rm div}_X(a) = \pi_*({\rm div}_Z(a)) = 
nq_1\pi_*([W_1]) + \cdots + nq_l\pi_*([W_l]) 
\]
by Proposition~1.4 in \cite{F}.
Since $\pi(W_i) \subset Y$, $\codim_X\pi(W_i)$ is at least $2$ for each $i$.
Therefore, we have ${\rm div}_X(a) = 0$.
Since $X$ is an affine normal scheme, 
$a$ is a unit of $\Gamma(X, {\cal O}_X)$.
Hence $a$ is also a unit of $\Gamma(Z, {\cal O}_Z)$
and ${\rm div}_Z(a) = 0$.
It follows that $q_1 = \cdots = q_l = 0$.
\qed

\noindent
\begin{tabular}{l}
Department of Mathematics\\
Meiji University\\
Higashi-Mita 1-1-1, Tama-ku,\\
Kawasaki 214-8571, Japan\\
{\tt kurano@math.meiji.ac.jp} \\
{\tt http://www.math.meiji.ac.jp/\~{}kurano}
\end{tabular}

\end{document}